\documentclass[a4paper,12pt]{article}
\usepackage[cp1251]{inputenc} 
\usepackage[english]{babel}
\usepackage{amsfonts,amssymb}
\title{The homotopy invariance of dihedral homology of involutive $\A$-algebras over rings.}
\author{S.V. Lapin}
\date{}
\newcommand{\F}{F_{\infty}}
\newcommand{\D}{D_{\infty}}
\newcommand{\bu}{\bullet}
\newcommand{\E}{E_{\infty}}
\newcommand{\A}{A_{\infty}}

\newcommand{\p}{\partial}

\begin{document}
\maketitle

\begin{abstract} The dihedral homology
functor $HD:A_\infty^{{\rm inv}}(K)\to GrM(K)$ from the category
$A_\infty^{{\rm inv}}(K)$ of involutive $A_\infty$-al\-geb\-ras
over any commutative unital ring $K$ to the category $GrM(K)$ of
graded $K$-modules is constructed. Further, it is showed that this
functor sends homotopy equivalences of involutive
$A_\infty$-al\-geb\-ras into isomorphisms of graded modules.
\end{abstract}
\vspace{0.5cm}

In \cite{SLapin1}, on the basis of the combinatorial and homotopy
technique of differential modules with $\infty$-simplicial faces
\cite{Lap1}-\cite{Lap7} and $\D$-differential modules
\cite{Lap9}-\cite{Lap17} the dihedral triple complex of an
involutive $\A$-algebra over any commutative unital ring was
constructed. This triple complex generalizes the dihedral triple
complex \cite{KLS1}-\cite{Sol} of an involutive associative
algebra given over an arbitrary commutative unital ring. Further,
in \cite{SLapin1}, dihedral homology of any involutive
$\A$-al\-geb\-ra over an arbitrary commutative unital ring was
defined as the homology of the chain complex associated with the
dihedral triple complex of this involutive $\A$-algebra. The
dihedral homology of involutive $\A$-algebras over commutative
unital rings introduced in \cite{SLapin1} generalizes the dihedral
homology of involutive associative algebras over commutative
unital rings defined in \cite{KLS1}. It is well known
\cite{KLS1}-\cite{Sol} that over fields of characteristic zero the
dihedral homology introduced in \cite{KLS1} is isomorphic to the
dihedral homology defined in \cite{Tsyg} by using the complex of
coinvariants for the action of dihedral groups. Similar to this,
in \cite{SLapin1}, it was shown that over fields of characteristic
zero the dihedral homology of involutive $\A$-algebras introduced
in \cite{SLapin1} is isomorphic to the dihedral homology of
involutive $\A$-algebras defined in \cite{B} by using the complex
of coinvariants for the action of dihedral groups. Moreover, in
\cite{SLapin1}, for involutive homotopy unital $\A$-algebras over
any commutative unital rings, the analogue of the
Krasauskas-Lapin-Solov'ev exact sequence was constructed.

On the other hand, in \cite{Kad}, it was established that the
homology of the $B$-con\-struc\-tion of an $\A$-al\-geb\-ra is
homotopy invariant, i.e., this homology is invariant under
homotopy equivalences of $\A$-al\-geb\-ras. Now note that the
dihedral triple complex of an involutive $\A$-al\-geb\-ra
constructed in \cite{SLapin1} is the dihedral analogue of the
$B$-con\-struc\-tion of an $\A$-al\-geb\-ra, and the dihedral
homology of an involutive $\A$-al\-geb\-ra is defined in
\cite{SLapin1} as the homology of this dihedral analogue of the
$B$-con\-struc\-tion. This gives rise to an interesting natural
question: do the dihedral homology of involutive $\A$-al\-geb\-ras
is homotopy invariant under the homotopy equivalences of
involutive $\A$-al\-geb\-ras? In present paper a positive answer
to this question is given.

The paper consists of three paragraphs. In the first paragraph, we
first recall necessary definitions related to the notion of a
dihedral module with $\infty$-sim\-pli\-cial faces or, more
briefly, an $D\F$-module \cite{SLapin1}, which homotopy
generalizes the notion of a dihedral module with simplicial faces
\cite{KLS1}-\cite{Sol}. After that, the category of
$D\F$-mo\-du\-les is defined, namely, the notion of a morphism of
$D\F$-modules is introduced, and it is shown that the composition
of morphisms of $D\F$-mo\-du\-les is a morphism of
$D\F$-mo\-du\-les. Next, the concept of a homotopy between
morphisms of $D\F$-mo\-du\-les and the notion of a homotopy
equivalence of $D\F$-mo\-du\-les are introduced.

In the second paragraph, we first recall necessary definitions
related to the notion of a dihedral homology of $D\F$-modules
\cite{SLapin1}. Next, it is shown that the dihedral homology of
$D\F$-modules defines the functor from the category of
$D\F$-modules to the category of graded modules. In addition, it
is shown that this functor sends homotopy equivalences of
$D\F$-modules into isomorphisms of graded modules.

In the third paragraph, we first recall necessary definitions
related to the notion of an $\A$-algebra \cite{Kad} and an
involutive $\A$-algebra \cite{B}. After that we give the
definitions of a morphism of involutive $\A$-algebras and a
homotopy between morphisms of involutive $\A$-algebras and,
moreover, we introduce the notion of a homotopy equivalence of
involutive $\A$-algebras. Next, we recall the concept of a
dihedral homology of involutive $\A$-algebras over an arbitrary
commutative unital rings \cite{SLapin1}. Then, by using results of
the second paragraph, it is shown that the dihedral homology of
involutive $\A$-algebras defines the functor from the category of
involutive $\A$-algebras to the category of graded modules.
Moreover, it is shown that this functor sends homotopy
equivalences of involutive $\A$-algebras into isomorphisms of
graded modules. As a corollary, we obtain that the dihedral
homology of an involutive $\A$-algebra over any field is
isomorphic to the dihedral homology of the involutive $\A$-algebra
of homologies for the source involutive $\A$-algebra. In
particular, it is obtained that the dihedral homology of an
involutive associative differential algebra over any field is
isomorphic to the dihedral homology of the involutive $\A$-algebra
of homologies for the source involutive associative differential
algebra.

We proceed to precise definitions and statements. All modules and
maps of modules considered in this paper are, respectively,
$K$-modules and $K$-linear maps of modules, where $K$ is any
unital (i.e., with unit) commutative ring.
\vspace{0.5cm}

\centerline{\bf \S\,1. Dihedral modules with $\infty$-simplicial
faces and} \centerline{\bf their morphisms and homotopies}
\vspace{0.5cm}

In what follows, by a bigraded module we mean any bigraded module
$X=\{X_{n,\,m}\}$, $n\geqslant 0$, $m\geqslant 0$, and by a
differential bigraded module, or, briefly, a differential module
$(X,d)$, we mean any bigraded module $X$ endowed with a
differential $d:X_{*,\bu}\to X_{*,\bu-1}$ of bidegree $(0,-1)$.

Recall that a differential module with simplicial faces is defined
as a differential module $(X,d)$ together with a family of module
maps $\p_i:X_{n,\bu}\to X_{n-1,\bu}$, $0\leqslant i\leqslant n$,
which are maps of differential modules and satisfy the simplicial
commutation relations $\p_i\p_j=\p_{j-1}\p_i$, $i<j$. The maps
$\p_i:X_{n,\bu}\to X_{n-1,\bu}$ are called the simplicial face
operators or, more briefly, the simplicial faces of the
differential module $(X,d)$.

Now, we recall the notion of a differential module with
$\infty$-simplicial faces \cite{Lap1} (see also
\cite{Lap2}-\cite{Lap7}), which is a homotopy invariant analogue
of the notion of a differential module with simplicial faces.

Let $\Sigma_k$ be the symmetric group of permutations on a
$k$-element set. Given an arbitrary permutation
$\sigma\in\Sigma_k$ and any $k$-tuple of nonnegative integers
$(i_1,\dots,i_k)$, where $i_1<\dots<i_k$, we consider the
$k$-tuple $(\sigma(i_1),\dots,\sigma(i_k))$, where $\sigma$ acts
on the $k$-tuple $(i_1,\dots,i_k)$ in the standard way, i.e.,
permutes its components. For the $k$-tuple
$(\sigma(i_1),\dots,\sigma(i_k))$, we define a $k$-tuple
$(\widehat{\sigma(i_1)},\dots,\widehat{\sigma(i_k)})$ by the
following formulae $$\widehat{\sigma
(i_s)}=\sigma(i_s)-\alpha(\sigma(i_s)),\quad 1\leqslant s\leqslant
k,$$ where each $\alpha(\sigma(i_s))$ is the number of those
elements of $(\sigma(i_1),\dots,\sigma(i_s),\dots\sigma(i_k))$ on
the right of $\sigma(i_s)$ that are smaller than $\sigma(i_s)$.

A differential module with $\infty$-simplicial faces or, more
briefly, an $\F$-module $(X,d,\p)$ is defined as a differential
module $(X,d)$ together with a family of module maps
$$\p=\{\p_{(i_1,\dots ,i_k)}:X_{n,\bu}\to X_{n-k,\bu+k-1}\},\quad
1\leqslant k\leqslant n,$$
$$i_1,\dots,i_k\in\mathbb{Z},\quad
0\leqslant i_1<\dots<i_k\leqslant n,$$ which satisfy the relations
$$d(\p_{(i_1,\dots,i_k)})=\sum_{\sigma\in\Sigma_k}\sum_{I_{\sigma}}
(-1)^{{\rm sign}(\sigma)+1}
\p_{(\widehat{\sigma(i_1)},\dots,\widehat{\sigma(i_m)})}\,
\p_{(\widehat{\sigma(i_{m+1})},\dots,\widehat{\sigma(i_k)})},\eqno(1.1)$$
where $I_\sigma$ is the set of all partitions of the $k$-tuple
$(\widehat{\sigma(i_1)},\dots,\widehat{\sigma(i_k)})$ into two
tuples $(\widehat{\sigma(i_1)},\dots,\widehat{\sigma(i_m)})$ and
$(\widehat{\sigma(i_{m+1})},\dots,\widehat{\sigma(i_k)})$,
$1\leqslant m\leqslant k-1$, such that the conditions
$\widehat{\sigma(i_1)}<\dots<\widehat{\sigma(i_m)}$ and
$\widehat{\sigma(i_{m+1})}<\dots<\widehat{\sigma(i_k)}$ holds.

The family of maps $\p=\{\p_{(i_1,\dots ,i_k)}\}$ is called the
$\F$-differential of the $\F$-module $(X,d,\widetilde{\p})$. The
maps $\p_{(i_1,\dots ,i_k)}$ that form the $\F$-differential of an
$\F$-module $(X,d,\p)$ are called the $\infty$-simplicial faces of
this $\F$-module.

It is easy to show that, for $k=1,2,3$, relations $(1.1)$ take,
respectively, the following view $$d(\p_{(i)})=0,\quad i\geqslant
0,\quad d(\p_{(i,j)})=\p_{(j-1)}\p_{(i)}-\p_{(i)}\p_{(j)},\quad
i<j,$$
$$d(\p_{(i_1,i_2,i_3)})=-\p_{(i_1)}\p_{(i_2,i_3)}-\p_{(i_1,i_2)}\p_{(i_3)}-
\p_{(i_3-2)}\p_{(i_1,i_2)}-$$
$$-\,\p_{(i_2-1,i_3-1)}\p_{(i_1)}+\p_{(i_2-1)}\p_{(i_1,i_3)}+\p_{(i_1,i_3-1)}\p_{(i_2)},
\quad i_1<i_2<i_3.$$

It is easy to check that, for any permutation $\sigma\in\Sigma_k$
and any $k$-tuple $(i_1,\dots,i_k)$, where $i_1<\dots<i_k$, the
conditions $\widehat{\sigma(i_1)}<\dots<\widehat{\sigma(i_m)}$ and
$\widehat{\sigma(i_{m+1})}<\dots<\widehat{\sigma(i_k)}$ are
equivalent to the conditions $\sigma(i_1)<\dots<\sigma(i_m)$ and
$\sigma(i_{m+1})<\dots<\sigma(i_k)$. This readily implies that the
$k$-tuple
$(\widehat{\sigma(i_{m+1})},\dots,\widehat{\sigma(i_k)})$, which
specified in $(1.1)$, coincides with the $k$-tuple
$(\sigma(i_{m+1}),\dots,\sigma(i_k))$.

Simplest examples of differential modules with $\infty$-simplicial
faces are differential modules with simplicial faces. Indeed,
given any differential module with simplicial faces $(X,d,\p_i)$,
we can define the $\F$-differential $\p=\{\p_{(i_1,\dots ,i_k)}\}$
by setting $\p_{(i)}=\p_i$, $i\geqslant 0$, and
$\p_{(i_1,\dots,i_k)}=0$, $k>1$, thus obtaining the differential
module with $\infty$-sim\-pli\-cial faces $(X,d,\p)$.

It is worth mentioning that the notion of an differential module
with $\infty$-simplicial faces specified above is a part of the
general notion of a differential $\infty$-simplicial module
introduced in \cite{Lap3} by using the homotopy technique of
differential Lie modules over curved colored coalgebras.

Now we recall \cite{KLS1}-\cite{Sol} that a dihedral differential
module with simplicial faces $(X,d,\p_i,t,r)$ is defined as a
differential module with simplicial faces $(X,d,\p_i)$ equipped
with two families of module maps $$t=\{t_n:X_{n,\bu}\to
X_{n,\bu}\},\quad r=\{r_n:X_{n,\bu}\to X_{n,\bu}\},\quad
n\geqslant 0,$$ which satisfy the following relations:
$$t_n^{\,n+1}=1_{X_{n,\bu}},\quad r_n^2=1_{X_{n,\bu}},\quad
r_nt_n=t_n^{-1}r_n,$$ $$dt_n=t_nd,\quad
\p_it_n=\left\{\begin{array}{ll} t_{n-1}\p_{i-1},& i>0,\\
\p_n,&i=0,\\
\end{array}\right.$$
$$dr_n=r_nd,\quad \p_ir_n=r_{n-1}\p_{n-i},\quad 0\leqslant
i\leqslant n.$$

Note that if in the definition of a dihedral module with
simplicial faces we remove the family of automorphisms
$r_n:X_{n,\bu}\to X_{n,\bu}$, $n\geqslant 0$, then we obtain the
definition of a cyclic module with simplicial faces \cite{Con}.

Now, let us recall \cite{SLapin1} that a dihedral differential
module with $\infty$-simplicial faces or, more briefly, an
$D\F$-module $(X,d,\p,t,r)$ is defined as any $\F$-module
$(X,d,\p)$ together with two families of module maps
$t=\{t_n:X_{n,\bu}\to X_{n,\bu}\}$, $n\geqslant 0$, and
$r=\{r_n:X_{n,\bu}\to X_{n,\bu}\}$, $n\geqslant 0$, which satisfy
the following relations: $$t_n^{\,n+1}=1_{X_{n,\bu}},\quad
r_n^2=1_{X_{n,\bu}},\quad r_nt_n=t_n^{-1}r_n,$$ $$dt_n=t_nd,\quad
\p_{(i_1,\dots,i_k)}t_n=\left\{\begin{array}{ll}
t_{n-k}\p_{(i_1-1,\dots,i_k-1)},& i_1>0,\\
(-1)^{k-1}\p_{(i_2-1,\dots,i_k-1,n)},& i_1=0,\\
\end{array}\right.\eqno(1.2)$$
$$dr_n=r_nd,\quad
\p_{(i_1,\dots,i_k)}r_n=(-1)^{k(k-1)/2}r_{n-k}\p_{(n-i_k,\dots,n-i_1)}\eqno(1.3)$$

Note that if in the definition of a $D\F$-module we remove the
family of automorphisms $r_n:X_{n,\bu}\to X_{n,\bu}$, $n\geqslant
0$, then we obtain the definition of a $C\F$-module, i.e., cyclic
module with $\infty$-simplicial faces \cite{Lapin} (see also
\cite{Lapin1}). Moreover, it is worth mentioning that the notion
of a $D\F$-module specified above is a part of the general notion
of a dihedral $\infty$-simplicial module introduced in
\cite{SLapin2} by using the homotopy technique of differential
modules over curved colored coalgebras.

The family of maps $\p=\{\p_{(i_1,\dots ,i_k)}\}$ is called the
$\F$-differential of the $D\F$-mo\-du\-le $(X,d,\p,t,r)$. The maps
$\p_{(i_1,\dots ,i_k)}$ are called the $\infty$-simplicial faces
of this $D\F$-module.

Simplest examples of $D\F$-modules are dihedral differential
modules with simplicial faces. Indeed, given any dihedral
differential module $(X,d,\p_i,t,r)$ with simplicial faces, we can
define the $\F$-differential $\p=\{\p_{(i_1,\dots ,i_k)}\}$ by
setting $\p_{(i)}=\p_i$, $i\geqslant 0$, and
$\p_{(i_1,\dots,i_k)}=0$, $k>1$, thus obtaining the $D\F$-module
$(X,d,\p,t,r)$.

Now, we recall that a map $f:(X,d,\p_i)\to (Y,d,\p_i)$ of
differential modules with simplicial faces is defined as a map of
differential modules $f:(X,d)\to(Y,d)$ that satisfies the
relations $\p_if=f\p_i$, $i\geqslant 0$.

Let us consider the notion of a morphism of differential modules
with $\infty$-simplicial faces \cite{Lap1} (see also \cite{Lap4}),
which homotopically generalizes the notion of a map differential
modules with simplicial faces.

A morphism of $\F$-modules $f:(X,d,\p)\to (Y,d,\p)$ is defined as
a family of module maps $$f=\{f_{(i_1,\dots,i_k)}:X_{n,\bu}\to
Y_{n-k,\bu+k}\},\quad 0\leqslant k\leqslant n,$$
$$i_1,\dots,i_k\in\mathbb{Z},\quad 0\leqslant
i_1<\dots<i_k\leqslant n,$$ (at $k=0$ we will use the denotation
$f_{(\,\,)}$), which satisfy the relations
$$d(f_{(i_1,\dots,i_k)})=-\p_{(i_1,\dots,i_k)}f_{(\,\,)}+f_{(\,\,)}\p_{(i_1,\dots,i_k)}\,+$$
$$+\sum_{\sigma\in\Sigma_k}\sum_{I_{\sigma}}(-1)^{{\rm
sign}(\sigma)+1}\p_{(\widehat{\sigma(i_1)},\dots,
\widehat{\sigma(i_m)})}f_{(\widehat{\sigma(i_{m+1})},\dots,
\widehat{\sigma(i_k)})}\,-$$ $$-\,f_{(\widehat{\sigma(i_1)},\dots,
\widehat{\sigma(i_m)})}\p_{(\widehat{\sigma(i_{m+1})},\dots,
\widehat{\sigma(i_k)})},\eqno(1.4)$$ where $I_{\sigma}$ is the
same as in $(1.1)$. The maps $f_{(i_1,\dots ,i_k)}\in f$ are
called the components of the morphism $f:(X,d,\p)\to (Y,d,\p)$.

For example, at $k=0,1,2,3$ the relations $(1.4)$ take,
respectively, the following view $$d(f_{(\,\,)})=0,\qquad
d(f_{(i)})=f_{(\,\,)}\p_{(i)}-\p_{(i)}f_{(\,\,)},\quad i\geqslant
0,$$
$$d(f_{(i,j)})=-\p_{(i,j)}f_{(\,\,)}+f_{(\,\,)}\p_{(i,j)}-\p_{(i)}f_{(j)}+
\p_{(j-1)}f_{(i)}+ f_{(i)}\p_{(j)}-f_{(j-1)}\p_{(i)},\quad i<j,$$
$$d(f_{(i_1,i_2,i_3)})=-\p_{(i_1,i_2,i_3)}f_{(\,\,)}+f_{(\,\,)}\p_{(i_1,i_2,i_3)}-
\p_{(i_1)}f_{(i_2,i_3)}-\p_{(i_1,i_2)}f_{(i_3)}-
\p_{(i_3-2)}f_{(i_1,i_2)}-$$
$$-\,\p_{(i_2-1,i_3-1)}f_{(i_1)}+\p_{(i_2-1)}f_{(i_1,i_3)}+\p_{(i_1,i_3-1)}f_{(i_2)}+
f_{(i_1)}\p_{(i_2,i_3)}+f_{(i_1,i_2)}\p_{(i_3)}+$$
$$+f_{(i_3-2)}\p_{(i_1,i_2)}+
f_{(i_2-1,i_3-1)}\p_{(i_1)}-f_{(i_2-1)}\p_{(i_1,i_3)}-f_{(i_1,i_3-1)}\p_{(i_2)},\quad
i_1<i_2<i_3.$$

Now, we recall \cite{Lap1} that a composition of an arbitrary
given morphisms of $\F$-mo\-du\-les $f:(X,d,\p)\to (Y,d,\p)$ and
$g:(Y,d,\p)\to (Z,d,\p)$ is defined as  a morphism of $\F$-modules
$gf:(X,d,\p)\to (Z,d,\p)$ whose components are defined by
$$(gf)_{(i_1,\dots,i_k)}=\sum_{\sigma\in\Sigma_k}\sum_{I'_{\sigma}}
(-1)^{{\rm sign}(\sigma)}
g_{(\widehat{\sigma(i_1)},\dots,\widehat{\sigma(i_m)})}
f_{(\widehat{\sigma(i_{m+1})},\dots,\widehat{\sigma(i_k)})},\eqno(1.5)$$
where $I'_\sigma$ is the set of all partitions of the $k$-tuple
$(\widehat{\sigma(i_1)},\dots,\widehat{\sigma(i_k)})$ into two
tuples $(\widehat{\sigma(i_1)},\dots,\widehat{\sigma(i_m)})$ and
$(\widehat{\sigma(i_{m+1})},\dots,\widehat{\sigma(i_k)})$,
$0\leqslant m\leqslant k$, such that the conditions
$\widehat{\sigma(i_1)}<\dots<\widehat{\sigma(i_m)}$ and
$\widehat{\sigma(i_{m+1})}<\dots<\widehat{\sigma(i_k)}$ holds.

For example, at $k=0,1,2,3$ the formulae $(1.5)$ take,
respectively, the following form
$$(gf)_{(\,\,)}=g_{(\,\,)}f_{(\,\,)},\qquad(gf)_{(i)}=g_{(\,\,)}f_{(i)}+g_{(i)}f_{(\,\,)},$$
$$(gf)_{(i_1,i_2)}=g_{(\,\,)}f_{(i_1,i_2)}+g_{(i_1,i_2)}f_{(\,\,)}+
g_{(i_1)}f_{(i_2)}-g_{(i_2-1)}f_{(i_1)},\qquad i_1<i_2,$$
$$(gf)_{(i_1,i_2,i_3)}=g_{(\,\,)}f_{(i_1,i_2,i_3)}+g_{(i_1,i_2,i_3)}f_{(\,\,)}+
g_{(i_1)}f_{(i_2,i_3)}+g_{(i_1,i_2)}f_{(i_3)}+$$
$$+\,g_{(i_3-2)}f_{(i_1,i_2)}+
g_{(i_2-1,i_3-1)}f_{(i_1)}-g_{(i_2-1)}f_{(i_1,i_3)}-g_{(i_1,i_3-1)}f_{(i_2)},\quad
i_1<i_2<i_3.$$

{\bf Definition 1.1}. We define a morphism of $D\F$-modules
$$f:(X,d,\p,t,r)\to (Y,d,\p,t,r)$$ as any morphism of $\F$-modules
$f:(X,d,\p)\to (Y,d,\p)$ whose components satisfy the following
conditions: $$f_{(\,\,)}t_n=t_nf_{(\,\,)},\quad
f_{(i_1,\dots,i_k)}t_n=\left\{\begin{array}{lll}
t_{n-k}f_{(i_1-1,\dots,i_k-1)},&k\geqslant 1,& i_1>0,\\
(-1)^{k-1}f_{(i_2-1,\dots,i_k-1,n)},&k\geqslant 1,& i_1=0,\\
\end{array}\right.\eqno(1.6)$$
$$f_{(\,\,)}r_n=r_nf_{(\,\,)},\quad
f_{(i_1,\dots,i_k)}r_n=(-1)^{k(k-1)/2}r_{n-k}f_{(n-i_k,\dots,n-i_1)}\eqno(1.7)$$

Note that if in the definition of a morphism of $D\F$-modules we
remove the family of automorphisms $r_n:X_{n,\bu}\to X_{n,\bu}$,
$n\geqslant 0$, then we obtain the definition of a morphism of
$C\F$-modules, i.e., a morphism of cyclic modules with
$\infty$-simplicial faces \cite{Lapin2}. It means that any
morphism of $D\F$-modules $f:(X,d,\p,t,r)\to (Y,d,\p,t,r)$ always
is the morphism of $C\F$-modules $f:(X,d,\p,t)\to (Y,d,\p,t)$.

By using the fact that any morphism of $D\F$-modules is a morphism
of $\F$-mo\-du\-les we define the composition of morphisms of
$D\F$-modules as a composition of morphisms of $\F$-modules.

{\bf Theorem 1.1}. The composition of morphisms of $D\F$-modules
is a morphism of $D\F$-modules.

{\bf Proof}. Given any morphisms of $D\F$-modules
$f:(X,d,\p,t,r)\to (Y,d,\p,t,r)$ and $g:(Y,d,\p,t,r)\to
(Z,d,\p,t,r)$, we need to check that components of the morphism of
$\F$-modules $gf:(X,d,\p)\to (Y,d,\p)$ satisfy the relations
$(1.5)$ and $(1.6)$. In \cite{Lapin2} it was shown that the
composition of morphisms of $C\F$-modules is a morphism of
$C\F$-modules. It implies that components of the morphism of
$\F$-modules $gf:(X,d,\p)\to (Y,d,\p)$ satisfy the relations
$(1.5)$. Now, we check that the module maps
$(gf)_{(i_1,\dots,i_k)}\in gf$ satisfy the relations $(1.6)$. It
is clearly that at $k=0$ we have
$(gf)_{(\,\,)}r_n=r_n(gf)_{(\,\,)}$. At $k\geqslant 1$ by the
definition of a composition of morphisms of $\F$-modules we have
$$(gf)_{(i_1,\dots,i_k)}r_n=\sum_{\sigma\in\Sigma_k}\sum_{I'_{\sigma}}
(-1)^{{\rm sign}(\sigma)}
g_{(\widehat{\sigma(i_1)},\dots,\widehat{\sigma(i_m)})}
f_{(\widehat{\sigma(i_{m+1})},\dots,\widehat{\sigma(i_k)})}r_n,\eqno(1.8)$$
$$(-1)^{k(k-1)/2}r_{n-k}(gf)_{(n-i_k,\dots,n-i_1)}=$$
$$=(-1)^{k(k-1)/2}r_{n-k}\sum_{\varrho\in\Sigma_k}\sum_{I'_{\varrho}}
(-1)^{{\rm sign}(\varrho)}
g_{(\widehat{\varrho(n-i_k)},\dots,\widehat{\varrho(n-i_{k-m+1})})}
f_{(\widehat{\varrho(n-i_{k-m})},\dots,\widehat{\varrho(n-i_1)})}.\eqno(1.9)$$
Let us show that each summand on the right-hand side of $(1.8)$ is
equal to some summand on the right-hand side of $(1.9)$. Given any
fixed permutation $\sigma\in\Sigma_k$, consider the summand
$$(-1)^{{\rm sign}(\sigma)}
g_{(\widehat{\sigma(i_1)},\dots,\widehat{\sigma(i_m)})}
f_{(\widehat{\sigma(i_{m+1})},\dots,\widehat{\sigma(i_k)})}r_n$$
on the right-hand side of $(1.8)$. By using the relations $(1.7)$
we obtain $$(-1)^{{\rm sign}(\sigma)}
g_{(\widehat{\sigma(i_1)},\dots,\widehat{\sigma(i_m)})}
f_{(\widehat{\sigma(i_{m+1})},\dots,\widehat{\sigma(i_k)})}r_n=$$
$$=(-1)^{{\rm
sign}(\sigma)+\varepsilon}r_{n-k}g_{(n-(k-m)-\widehat{\sigma(i_m)},\dots,n-(k-m)-\widehat{\sigma(i_1)})}
f_{(n-\widehat{\sigma(i_{k})},\dots,n-\widehat{\sigma(i_{m+1})})},$$
where $$\varepsilon=\frac{m(m-1)}{2}+\frac{(k-m)(k-m-1)}{2}\,.$$
Now, given the permutation $\sigma\in\Sigma_k$ and the partition
$$(\widehat{\sigma(i_1)},\dots,\widehat{\sigma(i_m)}\,|\,\widehat{\sigma(i_{m+1})},\dots,\widehat{\sigma(i_k)})
\in I'_{\sigma},$$ we define the permutation $\varrho\in\Sigma_k$
of the collection $(n-i_k,n-i_{k-1},\dots,n-i_2,n-i_1)$ by the
following formulae:
$$\varrho(n-i_k)=n-\sigma(i_m),\varrho(n-i_{k-1})=n-\sigma(i_{m-1}),\dots,\varrho(n-i_{k-m+1})=n-\sigma(i_1),$$
$$\varrho(n-i_{k-m})=n-\sigma(i_k),\varrho(n-i_{k-m-1})=n-\sigma(i_{k-1}),\dots,
\varrho(n-i_1)=n-\sigma(i_{m+1}).$$ For more clarity, it is worth
saying that the permutation $\varrho$ is the product
$\varrho=\gamma\psi\vartheta$ of the following permutations:
$$(n-i_k,n-i_{k-1},\dots,n-i_2,n-i_1)\buildrel{\vartheta}\over\longrightarrow
(n-i_1,n-i_2,\dots,n-i_{k-1},n-i_k)\to$$
$$\buildrel{\psi}\over\longrightarrow
(n-\sigma(i_1),n-\sigma(i_2),\dots,n-\sigma(i_{k-1}),n-\sigma(i_k))=$$
$$=(n-\sigma(i_1),n-\sigma(i_2),\dots,,n-\sigma(i_m)\,|\,n-\sigma(i_{m+1}),n-\sigma(i_{m+2}),\dots,n-\sigma(i_k))\to$$
$$\buildrel{\gamma}\over\longrightarrow
(n-\sigma(i_m),n-\sigma(i_{m-1}),\dots,n-\sigma(i_1)\,|\,n-\sigma(i_k),n-\sigma(i_{k-1}),\dots,n-\sigma(i_{m+1})).$$
Comparing the tuples $(\sigma(i_1),\sigma(i_2),\dots,\sigma(i_k))$
and $(\varrho(n-i_k),\varrho(n-i_{k-1}),\dots,\varrho(n-i_1))$, we
see that
$$\widehat{\varrho(n-i_k)}=n-(k-m)-\widehat{\sigma(i_m)},\dots,\widehat{\varrho(n-i_{k-m+1})}=
n-(k-m)-\widehat{\sigma(i_1)},$$
$$\widehat{\varrho(n-i_{k-m})}=n-\widehat{\sigma(i_k)},\dots,\widehat{\varrho(n-i_1)}=
n-\widehat{\sigma(i_{m+1})},$$ $${\rm sign}(\varrho)={\rm
sign}(\gamma\psi\vartheta)={\rm sign}(\gamma)+{\rm
sign}(\psi)+{\rm sign}(\vartheta)=$$ $$=\frac{k(k-1)}{2}+{\rm
sign}(\sigma)+\frac{m(m-1)}{2}+\frac{(k-m)(k-m-1)}{2}\,.$$ Since
$(\widehat{\varrho(n-i_k)},\dots,\widehat{\varrho(n-i_{k-m+1})}\,|\,\widehat{\varrho(n-i_{k-m})},
\dots,\widehat{\varrho(n-i_1)})\in I'_{\varrho}$, the right-hand
side of $(1.9)$ contains the summand
$$(-1)^{k(k-1)/2}r_{n-k}(-1)^{{\rm sign}(\varrho)}
g_{(\widehat{\varrho(n-i_k)},\dots,\widehat{\varrho(n-i_{k-m+1})})}
f_{(\widehat{\varrho(n-i_{k-m})},\dots,\widehat{\varrho(n-i_1)})}.$$
It is clear that the following equality holds:
$$(-1)^{k(k-1)/2}r_{n-k}(-1)^{{\rm sign}(\varrho)}
g_{(\widehat{\varrho(n-i_k)},\dots,\widehat{\varrho(n-i_{k-m+1})})}
f_{(\widehat{\varrho(n-i_{k-m})},\dots,\widehat{\varrho(n-i_1)})}=$$
$$=(-1)^{{\rm
sign}(\sigma)+\varepsilon}r_{n-k}g_{(n-(k-m)-\widehat{\sigma(i_m)},\dots,n-(k-m)-\widehat{\sigma(i_1)})}
f_{(n-\widehat{\sigma(i_{k})},\dots,n-\widehat{\sigma(i_{m+1})})}.$$
Thus, we have shown that each summand on the right-hand side of
$(1.8)$ is equal to a summand on the right-hand side of $(1.9)$.
It follows that the right-hand sides of $(1.8)$ and $(1.9)$ are
equal, because the number of summands on the right-hand side of
$(1.8)$ equals that on the right-hand side of $(1.9)$ and,
moreover, the specified above permutations $\sigma$ and $\varrho$
uniquely determine one another.~~~$\blacksquare$

It is clear that the associativity of the composition of morphisms
of $\F$-modules implies the associativity of the composition of
morphisms of $D\F$-modules. Moreover, for each $D\F$-module
$(X,d,\p,t,r)$, there is the identity morphism
$$1_X=\{(1_X)_{(i_1,\dots,i_k)}\}:(X,d,\p,t,r)\to (X,d,\p,t,r),$$
where $(1_X)_{(\,\,)}={\rm id_X}$ and $(1_X)_{(i_1,\dots,i_k)}=0$
for all $k\geqslant 1$. Thus, the class of all $D\F$-mo\-du\-les
over any commutative unital ring $K$ and their morphisms is a
category, which we denote by $D\F(K)$.

Now, we recall that a differential homotopy or, more briefly, a
homotopy between morphisms $f,g:(X,d,\p_i)\to(Y,d,\p_i)$ of
differential modules with simplicial faces is defines as a
differential homotopy $h:X_{*,\bu}\to Y_{*,\bu+1}$ between
morphisms of differential modules $f,g:(X,d)\to (Y,d)$, which
satisfies the relations $\p_ih+h\p_i=0$, $i\geqslant 0$.

Let us consider the notion of a homotopy between morphisms of
$\F$-modules \cite{Lap1} (see also \cite{Lap4}), which
homotopically generalizes the notion of a homotopy between
morphisms of differential modules with simplicial faces.

A homotopy between morphisms of $\F$-modules $f,g:(X,d,\p)\to
(Y,d,\p)$ is defined as a family of module maps
$$h=\{h_{(i_1,\dots,i_k)}:X_{n,\bu}\to Y_{n-k,\bu+k+1}\},\quad
0\leqslant k\leqslant n,$$ $$i_1,\dots,i_k\in\mathbb{Z},\quad
0\leqslant i_1<\dots<i_k\leqslant n,$$ (at $k=0$ we will use the
denotation $h_{(\,\,)}$), which satisfy the relations
$$d(h_{(i_1,\dots,i_k)})=f_{(i_1,\dots,i_k)}-g_{(i_1,\dots,i_k)}-
\partial_{(i_1,\dots,i_k)}h_{(\,\,)}-h_{(\,\,)}\partial_{(i_1,\dots,i_k)}\,+$$
$$+\sum_{\sigma\in\Sigma_k}\sum_{I_{\sigma}}(-1)^{{\rm
sign}(\sigma)+1}
\partial_{(\widehat{\sigma(i_1)},\dots,\widehat{\sigma(i_m)})}
h_{(\widehat{\sigma(i_{m+1})},\dots,\widehat{\sigma(i_k)})}\,+$$
$$+\,h_{(\widehat{\sigma(i_1)},\dots,\widehat{\sigma(i_m)})}
\partial_{(\widehat{\sigma(i_{m+1})},\dots,\widehat{\sigma(i_k)})},\eqno(1.10)$$
where $I_{\sigma}$ is the same as in $(1.1)$. The maps
$h_{(i_1,\dots ,i_k)}\in h$ are called the components of the
homotopy $h$.

For example, at $k=0,1,2,3$ the relations $(1.10)$ take,
respectively, the following view
$$d(h_{(\,\,)})=f_{(\,\,)}-g_{(\,\,)},\quad
d(h_{(i)})=f_{(i)}-g_{(i)}-\p_{(i)}h_{(\,\,)}-h_{(\,\,)}\p_{(i)},\quad
i\geqslant 0,$$ $$d(h_{(i,j)})=f_{(i,j)}-g_{(i,j)}-
\p_{(i,j)}h_{(\,\,)}-h_{(\,\,)}\p_{(i,j)}-\p_{(i)}h_{(j)}+\p_{(j-1)}h_{(i)}-
h_{(i)}\p_{(j)}+h_{(j-1)}\p_{(i)},~i<j,$$
$$d(h_{(i_1,i_2,i_3)})=f_{(i_1,i_2,i_3)}-g_{(i_1,i_2,i_3)}-\p_{(i_1,i_2,i_3)}h_{(\,\,)}-h_{(\,\,)}\p_{(i_1,i_2,i_3)}-
\p_{(i_1)}f_{(i_2,i_3)}-\p_{(i_1,i_2)}f_{(i_3)}\,-$$
$$-\,\p_{(i_3-2)}f_{(i_1,i_2)}-\p_{(i_2-1,i_3-1)}f_{(i_1)}+\p_{(i_2-1)}f_{(i_1,i_3)}+\p_{(i_1,i_3-1)}f_{(i_2)}-
h_{(i_1)}\p_{(i_2,i_3)}-h_{(i_1,i_2)}\p_{(i_3)}\,-$$
$$-\,h_{(i_3-2)}\p_{(i_1,i_2)}-
h_{(i_2-1,i_3-1)}\p_{(i_1)}+h_{(i_2-1)}\p_{(i_1,i_3)}+h_{(i_1,i_3-1)}\p_{(i_2)},\quad
i_1<i_2<i_3.$$

{\bf Definition 1.2}. We define a homotopy
$h=\{h_{(i_1,\dots,i_k)}\}$ between an arbitrary morphisms of
$D\F$-modules $f,g:(X,d,\p,t,r)\to (Y,d,\p,t,r)$ as any homotopy
$h=\{h_{(i_1,\dots,i_k)}\}$ between morphisms of $\F$-modules
$f,g:(X,d,\p)\to (Y,d,\p)$ whose components satisfy the following
conditions: $$h_{(\,\,)}t_n=t_nh_{(\,\,)},\quad
h_{(i_1,\dots,i_k)}t_n=\left\{\begin{array}{lll}
t_{n-k}h_{(i_1-1,\dots,i_k-1)},&k\geqslant 1,& i_1>0,\\
(-1)^{k-1}h_{(i_2-1,\dots,i_k-1,n)},&k\geqslant 1,& i_1=0,\\
\end{array}\right.\eqno(1.11)$$
$$h_{(\,\,)}r_n=r_nh_{(\,\,)},\quad
h_{(i_1,\dots,i_k)}r_n=(-1)^{k(k-1)/2}r_{n-k}h_{(n-i_k,\dots,n-i_1)}\eqno(1.12)$$

{\bf Proposition 1.1}. For any $D\F$-modules $(X,d,\p,t,r)$ and
$(Y,d,\p,t,r)$, the relation between morphisms of $D\F$-modules of
the form $(X,d,\p,t,r)\to (Y,d,\p,t,r)$ defined by the presence of
a homotopy between them is an equivalence relation.

{\bf Proof}. Suppose given an arbitrary morphism
$f:(X,d,\p,t,r)\to (Y,d,\p,t,r)$ of $D\F$-mo\-du\-les. Then we
have the homotopy $0=\{0_{(i_1,\dots,i_k)}=0\}$ between morphisms
of $D\F$-modules $f$ and $f$. Suppose given a homotopy
$h=\{h_{(i_1,\dots,i_k)}\}$ between morphisms of $D\F$-modules $f$
and $g$. Then the family of maps $-h=\{-h_{(i_1,\dots,i_k)}\}$ is
a homotopy between morphisms of $D\F$-modules $g$ and $f$. Suppose
given a homotopy $h=\{h_{(i_1,\dots,i_k)}\}$ between morphisms of
$D\F$-modules $f$ and $g$ and, moreover, given a homotopy
$H=\{H_{(i_1,\dots,i_k)}\}$ between morphisms of $D\F$-modules $g$
and $p$. Then the family of maps
$h+H=\{h_{(i_1,\dots,i_k)}+H_{(i_1,\dots,i_k)}\}$ is a homotopy
between morphisms of $D\F$-modules $f$ and $p$.~~~$\blacksquare$

By using specified in Proposition 1.1 the equivalence relation
between morphisms of $D\F$-modules the notion of a homotopy
equivalence of $D\F$-modules is introduced in the usual way.
Namely, a morphism of $D\F$-modules is called a homotopy
equivalence of $D\F$-modules, when this morphism have a homotopy
inverse morphism of $D\F$-mo\-dules.
\vspace{0.5cm}

\centerline{\bf \S\,2. The homotopy invariance of dihedral
homology of $D\F$-modules.}
\vspace{0.5cm}

First, recall that a $\D$-differential module \cite{Lap9} (sees
also \cite{Lap10}-\cite{Lap17}) or, more briefly, a $\D$-module
$(X,d^{\,i})$ is defined as a module $X$ together with a family of
module maps $\{d^{\,i}:X\to X~|~i\in\mathbb{Z},~i\geqslant 0\}$
satisfying the relations
$$\sum\limits_{i+j=k}d^{\,i}d^{\,j}=0,\quad k\geqslant
0.\eqno(2.1)$$

It is worth noting that a $\D$-module $(X,d^{\,i})$ can be
equipped with any $\mathbb{Z}^{\times n}$-gra\-ding, i.e.,
$X=\{X_{k_1,\dots,k_n}\}$, where
$(k_1,\dots,k_n)\in\mathbb{Z}^{\times n}$ and $n\geqslant 1$, and
the module maps $d^{\,i}:X\to X$ can have any $n$-degree
$(l_1(i),\dots,l_n(i))\in\mathbb{Z}^{\times n}$ for each
$i\geqslant 0$, i.e., $d^{\,i}:X_{k_1,\dots,k_n}\to
X_{k_1+l_1(i),\dots,k_n+l_n(i)}$.

For $k=0$, the relations $(2.1)$ have the form $d^{\,0}d^{\,0}=0$,
and hence $(X,d^{\,0})$ is a differential module. In \cite{Lap9}
the homotopy invariance of the $\D$-module structure over any
unital commutative ring under homotopy equivalences of
differential modules was established. Later, it was shown in
\cite{LV} that the homotopy invariance of the $\D$-mo\-du\-le
structure over fields of characteristic zero can be established by
using the Koszul duality theory.

It is also worth saying that in \cite{Lap9} by using specified
above homotopy invariance of the $\D$-differential module
structure the relationship between $\D$-differential modules and
spectral sequences was established. More precisely, in \cite{Lap9}
was shown that over an arbitrary field the category of
$\D$-differential modules is equivalent to the category of
spectral sequences.

Now, we recall \cite{Lap9} that a $\D$-module $(X,d^{\,i})$ is
said to be stable if, for any $x\in X$, there exists a number
$k=k(x)\geqslant 0$ such that $d^{\,i}(x)=0$ for each $i>k$. Any
stable $\D$-module $(X,d^{\,i})$ determines the differential
$\overline{d\,}:X\to X$ defined by
$\overline{d\,}=(d^{\,0}+d^{\,1}+\dots+d^{\,i}+\dots)$. The map
$\overline{d\,}:X\to X$ is indeed a differential because relations
$(2.1)$ imply the equality $\overline{d\,}\,\overline{d\,}=0$. It
is easy to see that if the stable $\D$-mo\-du\-le $(X,d^{\,i})$ is
equipped with a $\mathbb{Z}^{\times n}$-grading
$X=\{X_{k_1,\dots,k_n}\}$, where $k_1\geqslant
0,\dots,k_n\geqslant 0$, and maps $d^{\,i}:X\to X$, $i\geqslant
0$, have $n$-degree $(l_1(i),\dots,l_n(i))$ satisfying the
condition $l_1(i)+\dots+l_n(i)=-1$, then there is the chain
complex $(\overline{X},\overline{d\,})$ defined by the following
formulae:
$$\overline{X}_m=\bigoplus_{k_1+\dots+k_n=m}X_{k_1\dots,k_n},\quad
\overline{d\,}=\sum_{i=0}^\infty d^{\,i}:\overline{X}_m\to
\overline{X}_{m-1},\quad m\geqslant 0.$$

It was shown in \cite{Lap1} that any $\F$-module $(X,d,\p)$
determines the sequence of stable $\D$-modules $(X,d_q^{\,i})$,
$q\geqslant 0$, equipped with the bigrading $X=\{X_{n,m}\}$,
$n\geqslant 0$, $m\geqslant 0$, and defined by the following
formulae: $$d_q^{\,0}=d,~~ d_q^{\,k}=\sum\limits_{0\leqslant
i_1<\dots<i_k\leqslant n-q
}(-1)^{i_1+\dots+i_k}\p_{(i_1,\dots,i_k)}:X_{n,\bu}\to
X_{n-k,\bu+k-1},~~k\geqslant 1.\eqno(2.2)$$

Let us recall \cite{Lapin} the construction of the chain bicomplex
$(C(\overline{X}),\delta_1,\delta_2)$ that is defined by the
$C\F$-module $(X,d,\p,t)$. Given any $C\F$-mo\-dule $(X,d,\p,t)$,
consider the two $\D$-mo\-dules $(X,d_0^{\,i})$ and
$(X,d_1^{\,i})$ defined by $(2.2)$ for $q=0,1$, and the two
families of maps $$T_n=(-1)^nt_n:X_{n,\bu}\to X_{n,\bu},\quad
n\geqslant 0,$$ $$N_n=1+T_n+T_n^2+\dots+T_n^n:X_{n,\bu}\to
X_{n,\bu},\quad n\geqslant 0.$$ Moreover, we consider the chain
complexes $(\overline{X},b)$ and $(\overline{X},b^{'})$ that
corresponded to the specified above $\D$-mo\-dules $(X,d_0^{\,i})$
and $(X,d_1^{\,i})$, where
$$b=\overline{d\,}_0=(d_0^{\,0}+d_0^{\,1}+\dots
+d_0^{\,i}+\dots)\quad\mbox{and}\quad
b^{'}=\overline{d\,}_1=(d_1^{\,0}+d_1^{\,1}+\dots+d_1^{\,i}+\dots).$$
Also consider the two families of maps
$$\overline{T}_n=\sum_{k=0}^nT_{k}:\overline{X}_n\to\overline{X}_n,\quad
\overline{N}_n=\sum_{k=0}^nN_{k}:\overline{X}_n\to\overline{X}_n,\quad
n\geqslant 0.$$ In \cite{Lapin} it was shown that the following
relations holds: $$(1-\overline{T}_n)\overline{N}_n=0,\quad
\overline{N}_n(1-\overline{T}_n)=0,\quad n\geqslant 0,$$
$$b(1-\overline{T}_n)=(1-\overline{T}_{n-1})b^{'},\quad
b^{'}\overline{N}_n=\overline{N}_{n-1}b,\quad n\geqslant 0.$$ It
follows from these relations that any $C\F$-module $(X,d,\p,t)$
determines the chain bicomplex \vspace{1cm}

\begin{center}
\parbox {4cm}
{\setlength{\unitlength}{1cm}
\begin{picture}(2,3.5)

\put(-3.1,3.5){\makebox(0,0){\vspace{0.3cm}$\vdots$}}
\put(-0.7,3.5){\makebox(0,0){\vspace{0.3cm}$\vdots$}}
\put(1.7,3.5){\makebox(0,0){\vspace{0.3cm}$\vdots$}}
\put(4.1,3.5){\makebox(0,0){\vspace{0.3cm}$\vdots$}}

\put(-3.1,3.2){\vector(0,-1){0.7}}
\put(-0.7,3.2){\vector(0,-1){0.7}}
\put(1.7,3.2){\vector(0,-1){0.7}}
\put(4.1,3.2){\vector(0,-1){0.7}}

\put(-3,2){\makebox(0,0){$\overline{X}_{n+1}$}}
\put(-0.6,2){\makebox(0,0){$\overline{X}_{n+1}$}}
\put(1.8,2){\makebox(0,0){$\overline{X}_{n+1}$}}
\put(4.2,2){\makebox(0,0){$\overline{X}_{n+1}$}}
\put(6.3,2){\makebox(0,0){$\dots$}}

\put(-1.3,2){\vector(-1,0){1.1}} \put(1.1,2){\vector(-1,0){1.1}}
\put(3.5,2){\vector(-1,0){1.1}} \put(5.9,2){\vector(-1,0){1.1}}

\put(-3.3,2.8){\makebox(0,0){$^{b}$}}
\put(-3.3,1.1){\makebox(0,0){$^{b}$}}
\put(-3.3,-0.6){\makebox(0,0){$^{b}$}}
\put(-3.3,-2.3){\makebox(0,0){$^{b}$}}

\put(-1.0,2.8){\makebox(0,0){$^{-b^{'}}$}}
\put(-1.0,1.1){\makebox(0,0){$^{-b^{'}}$}}
\put(-1.0,-0.6){\makebox(0,0){$^{-b^{'}}$}}
\put(-1.0,-2.3){\makebox(0,0){$^{-b^{'}}$}}

\put(1.5,2.8){\makebox(0,0){$^{b}$}}
\put(1.5,1.1){\makebox(0,0){$^{b}$}}
\put(1.5,-0.6){\makebox(0,0){$^{b}$}}
\put(1.5,-2.3){\makebox(0,0){$^{b}$}}

\put(3.8,2.8){\makebox(0,0){$^{-b^{'}}$}}
\put(3.8,1.1){\makebox(0,0){$^{-b^{'}}$}}
\put(3.8,-0.6){\makebox(0,0){$^{-b^{'}}$}}
\put(3.8,-2.3){\makebox(0,0){$^{-b^{'}}$}}

\put(-3.1,1.7){\vector(0,-1){1.0}}
\put(-0.7,1.7){\vector(0,-1){1.0}}
\put(1.7,1.7){\vector(0,-1){1.0}}
\put(4.1,1.7){\vector(0,-1){1.0}}

\put(-3,0.2){\makebox(0,0){$\overline{X}_n$}}
\put(-0.6,0.2){\makebox(0,0){$\overline{X}_n$}}
\put(1.8,0.2){\makebox(0,0){$\overline{X}_n$}}
\put(4.2,0.2){\makebox(0,0){$\overline{X}_n$}}
\put(6.3,0.2){\makebox(0,0){$\dots$}}

\put(-1.3,0.2){\vector(-1,0){1.1}}
\put(1.1,0.2){\vector(-1,0){1.1}}
\put(3.5,0.2){\vector(-1,0){1.1}}
\put(5.9,0.2){\vector(-1,0){1.1}}

\put(-3.1,-0.1){\vector(0,-1){1.0}}
\put(-0.7,-0.1){\vector(0,-1){1.0}}
\put(1.7,-0.1){\vector(0,-1){1.0}}
\put(4.1,-0.1){\vector(0,-1){1.0}}

\put(-3,-1.55){\makebox(0,0){$\overline{X}_{n-1}$}}
\put(-0.6,-1.55){\makebox(0,0){$\overline{X}_{n-1}$}}
\put(1.8,-1.55){\makebox(0,0){$\overline{X}_{n-1}$}}
\put(4.2,-1.55){\makebox(0,0){$\overline{X}_{n-1}$}}
\put(6.3,-1.55){\makebox(0,0){$\dots$}}

\put(-1.7,-1.3){\makebox(0,0){$^{1-\overline{T}_{n-1}}$}}
\put(-1.7,0.4){\makebox(0,0){$^{1-\overline{T}_n}$}}
\put(-1.7,2.2){\makebox(0,0){$^{1-\overline{T}_{n+1}}$}}

\put(0.6,-1.3){\makebox(0,0){$^{\overline{N}_{n-1}}$}}
\put(0.6,0.4){\makebox(0,0){$^{\overline{N}_n}$}}
\put(0.6,2.2){\makebox(0,0){$^{\overline{N}_{n+1}}$}}

\put(3.1,-1.3){\makebox(0,0){$^{1-\overline{T}_{n-1}}$}}
\put(3.1,0.4){\makebox(0,0){$^{1-\overline{T}_n}$}}
\put(3.1,2.2){\makebox(0,0){$^{1-\overline{T}_{n+1}}$}}

\put(5.4,-1.3){\makebox(0,0){$^{\overline{N}_{n-1}}$}}
\put(5.4,0.4){\makebox(0,0){$^{\overline{N}_n}$}}
\put(5.4,2.2){\makebox(0,0){$^{\overline{N}_{n+1}}$}}

\put(-3.1,-2.8){\makebox(0,0){$\vdots$}}
\put(-0.7,-2.8){\makebox(0,0){$\vdots$}}
\put(1.7,-2.8){\makebox(0,0){$\vdots$}}
\put(4.1,-2.8){\makebox(0,0){$\vdots$}}

\put(-1.3,-1.55){\vector(-1,0){1.1}}
\put(1.1,-1.55){\vector(-1,0){1.1}}
\put(3.5,-1.55){\vector(-1,0){1.1}}
\put(5.9,-1.55){\vector(-1,0){1.1}}

\put(-3.1,-1.9){\vector(0,-1){0.7}}
\put(-0.7,-1.9){\vector(0,-1){0.7}}
\put(1.7,-1.9){\vector(0,-1){0.7}}
\put(4.1,-1.9){\vector(0,-1){0.7}}
\end{picture}}
\end{center}
\vspace{3.5cm}

\noindent We denote this chain bicomplex by
$(C(\overline{X}),\delta_1,\delta_2)$, where
$C(\overline{X})_{n,m}=\overline{X}_n$, $n\geqslant 0$,
$m\geqslant 0$, $\delta_1:C(\overline{X})_{n,m}\to
C(\overline{X})_{n-1,m}$, $\delta_2:C(\overline{X})_{n,m}\to
C(\overline{X})_{n,m-1}$, $$\delta_1=\left\{\begin{array}{ll}
b,&m\equiv 0\,{\rm mod}(2),\\ -b^{'},&m\equiv 1\,{\rm mod}(2),\\
\end{array}\right.\quad\delta_2=\left\{\begin{array}{ll}
1-\overline{T}_n,&m\equiv 1\,{\rm mod}(2),\\
\overline{N}_n,&m\equiv 0\,{\rm mod}(2).\\
\end{array}\right.$$
The chain complex associated with the chain bicomplex
$(C(\overline{X}),\delta_1,\delta_2)$ we denote by $({\rm
Tot}(C(\overline{X})),\delta)$, where $\delta=\delta_1+\delta_2$.

Recall \cite{Lapin} that the cyclic homology $HC(X)$ of a
$C\F$-module $(X,d,\p,t)$ is defined as the homology of the chain
complex $({\rm Tot}(C(\overline{X})),\delta)$ associated with the
chain bicomplex $(C(\overline{X}),\delta_1,\delta_2)$.

Now, suppose given any $D\F$-module $(X,d,\p,t,r)$. Since
$(X,d,\p,t,r)$ can be considered as the $C\F$-module $(X,d,\p,t)$,
we have the specified above chain complexes $(\overline{X},b)$ and
$(\overline{X},b^{'})$ and also the families of module maps
$\overline{T}_n:\overline{X}_n\to\overline{X}_n$, $n\geqslant 0$,
and $\overline{N}_n:\overline{X}_n\to\overline{X}_n$, $n\geqslant
0$. Given the $D\F$-module $(X,d,\p,t,r)$, consider the family of
module maps
$$\overline{R}_n=\sum_{m=0}^nR_m:\overline{X}_n\to\overline{X}_n,\quad
n\geqslant 0,$$ where $$R_m=(-1)^{m(m+1)/2}r_m:X_{m,\bu}\to
X_{m,\bu},\quad m\geqslant 0.$$ In \cite{SLapin1} it was shown
that the following relations holds:
$$\left.\begin{array}{c}(1-\overline{T}_n)(\overline{R}_n\overline{T}_n)=-\overline{R}_n(1-\overline{T}_n),
\quad\overline{N}_n\overline{R}_n=(\overline{R}_n\overline{T}_n)\overline{N}_n,\quad
n\geqslant 0,\\
\\
b\overline{R}_n=\overline{R}_{n-1}b,\quad
b^{'}(\overline{R}_n\overline{T}_n)=(\overline{R}_{n-1}\overline{T}_{n-1})b^{'},\quad
n\geqslant 0.\\
\end{array}\right\}\eqno(2.3)$$
Since the $D\F$-module $(X,d,\p,t,r)$ can be considered as the
$C\F$-module $(X,d,\p,t)$, the $D\F$-module $(X,d,\p,t,r)$ always
determines the chain bicomplex
$(C(\overline{X}),\delta_1,\delta_2)$. The relations $(2.3)$ say
us that there is a left action of the group
$\mathbb{Z}_2=\{1,\vartheta\,|\,\vartheta^2=1\}$ on the chain
bicomplex $(C(\overline{X}),\delta_1,\delta_2)$ of any
$D\F$-module $(X,d,\p,t,r)$. This left action is defined by means
of the automorphism $\vartheta:C(\overline{X})_{*,\bu}\to
C(\overline{X})_{*,\bu}$ of the order two, which, at any element
$x\in C(\overline{X})_{n,m}$, is given by the following rule:
$$\vartheta(x)=\left\{\begin{array}{ll}
(-1)^k\overline{R}_n(x),&m=2k,\\
(-1)^{k+1}\overline{R}_n\overline{T}_n(x),&m=2k+1.\\
\end{array}\right.$$

Recall \cite{SLapin1} that the dihedral homology $HD(X)$ of any
$D\F$-module $(X,d,\p,t,r)$ is defined as the hyperhomology
$H(\mathbb{Z}_2;(C(\overline{X}),\delta_1,\delta_2))$ of the group
$\mathbb{Z}_2$ with coefficients in the
$(C(\overline{X}),\delta_1,\delta_2)$ relative to the specified
above action of the group $\mathbb{Z}_2$ on the chain bicomplex
$(C(\overline{X}),\delta_1,\delta_2)$.

The hyperhomology
$H(\mathbb{Z}_2;(C(\overline{X}),\delta_1,\delta_2))$ is defined
as the homology of the chain complex that is associated with the
triple chain complex $({\cal
P}(\mathbb{Z}_2)\otimes_{K[\mathbb{Z}_2]}C(\overline{X}),\delta_1,\delta_2,\delta_3)$,
where $K[\mathbb{Z}_2]$ is a groups $K$-algebra of the group
$\mathbb{Z}_2=\{1,\vartheta\,|\,\vartheta^2=1\}$, the chain
complex $({\cal P}(\mathbb{Z}_2),d)$ is any projective resolvente
of the trivial $K[\mathbb{Z}_2]$-module $K$, and the differential
$\delta_3$ is defined by $$\delta_3=(-1)^{n+m}d\otimes 1:{\cal
P}(\mathbb{Z}_2)_l\otimes_{K[\mathbb{Z}_2]}C(\overline{X})_{n,m}\to
{\cal
P}(\mathbb{Z}_2)_{l-1}\otimes_{K[\mathbb{Z}_2]}C(\overline{X})_{n,m}.$$
If we take as the projective resolvente $({\cal
P}(\mathbb{Z}_2),d)$ the standard free resolvente $$({\cal
S}(\mathbb{Z}_2),d):\,K[\mathbb{Z}_2]\buildrel{~1-\vartheta}\over\longleftarrow
K[\mathbb{Z}_2]\buildrel{~1+\vartheta}\over\longleftarrow
K[\mathbb{Z}_2]\buildrel{~1-\vartheta}\over\longleftarrow
K[\mathbb{Z}_2]\buildrel{~1+\vartheta}\over\longleftarrow\cdots,$$
then we obtain that the dihedral homology
$HD(X)=H(\mathbb{Z}_2;(C(\overline{X}),\delta_1,\delta_2))$ of any
$D\F$-module $(X,d,\p,t,r)$ is the homology of the chain complex
that is associated with the triple chain complex
$(D(\overline{X}),\delta_1,\delta_2,\delta_3)$, where
$D(\overline{X})_{n,m,l}=C(\overline{X})_{n,m}=\overline{X}_n$,
$n\geqslant 0$, $m\geqslant 0$, $l\geqslant 0$, the differentials
$$\delta_1:D(\overline{X})_{n,m,l}=C(\overline{X})_{n,m}\to
C(\overline{X})_{n-1,m}=D(\overline{X})_{n-1,m,l},$$
$$\delta_2:D(\overline{X})_{n,m,l}=C(\overline{X})_{n,m}\to
C(\overline{X})_{n,m-1}=D(\overline{X})_{n,m-1,l}$$ were defined
above, and the differential $\delta_3:D(\overline{X})_{n,m,l}\to
D(\overline{X})_{n,m,l-1}$ is given by
$$\delta_3=\left\{\begin{array}{ll}
(-1)^n(1+(-1)^l\overline{R}_n),&m\equiv 0\,{\rm mod}(4),\\
(-1)^{n+1}(1+(-1)^{l+1}\overline{R}_n\overline{T}_n),&m\equiv
1\,{\rm mod}(4),\\ (-1)^n(1+(-1)^{l+1}\overline{R}_n),&m\equiv
2\,{\rm mod}(4),\\
(-1)^{n+1}(1+(-1)^l\overline{R}_n\overline{T}_n),&m\equiv 3\,{\rm
mod}(4).\\
\end{array}\right.$$
We denote by $({\rm Tot}(D(\overline{X})),\widehat{\delta})$,
where $\widehat{\delta}=\delta_1+\delta_2+\delta_3$, the chain
complex associated to the triple chain complex
$(D(\overline{X}),\delta_1,\delta_2,\delta_3)$.

Now, we investigate functorial and homotopy properties of the
dihedral homology of $D\F$-mo\-du\-les.

{\bf Theorem 2.1}. The dihedral homology of $D\F$-modules over any
commutative unital ring $K$ determines the functor $HD:D\F(K)\to
GrM(K)$ from the category of $D\F$-mo\-dules $D\F(K)$ to the
category of graded $K$-modules $GrM(K)$. This functor sends
homotopy equivalences of $D\F$-modules into isomorphisms of graded
modules.

{\bf Proof}. First, show that every morphism of $D\F$-modules
induces a module map of dihedral homology of these $D\F$-modules.
Given an arbitrary morphism of $D\F$-mo\-du\-les
$f=\{f_{(i_1,\dots,i_k}\}:(X,d,\p,t,r)\to (X,d,\p,t,r)$, consider
the family of module maps $$f_q^k=\sum\limits_{0\leqslant
i_1<\dots<i_k\leqslant n-q
}(-1)^{i_1+\dots+i_k}f_{(i_1,\dots,i_k)}:X_{n,\bu}\to
Y_{n-k,\bu+k},\quad k\geqslant 0,\quad q\geqslant 0.\eqno(2.4)$$
For the family of maps $\{f_q^k\}$, by using $(1.4)$ we obtain the
relations
$$\sum_{i+j=k}d_q^{\,i}f_q^j=\sum_{i+j=k}f_q^id_q^{\,j},\quad
k\geqslant 0,\quad q\geqslant 0,\eqno(2.5)$$ where $(X,d_q^{\,i})$
and $(Y,d_q^{\,i})$ are sequences of $\D$-modules respectively
defined by $(2.2)$ for $\F$-modules $(X,d,\p)$ and $(Y,d,\p)$.
Direct calculations with using $(1.6)$ show that the families maps
$\{f_0^k\}$ and $\{f_1^k\}$ satisfy the relations
$$f_0^k(1-T_n)=(1-T_{n-k})f_1^k,\quad f_1^kN_n=N_{n-k}f_0^k,\quad
k\geqslant 0,\quad n\geqslant 0.\eqno(2.6)$$ For $q=0,1$, the
equalities $(2.5)$ imply that maps of graded modules
$$\overline{f}_0=\sum_{k=0}^nf_0^k:\overline{X}_n\to\overline{Y}_n,\quad
\overline{f}_1=\sum_{k=0}^nf_1^k:\overline{X}_n\to\overline{Y}_n,\quad
n\geqslant 0,$$ are chain maps $\overline{f}_0:(\overline{X},b)\to
(\overline{Y},b)$, $\overline{f}_1:(\overline{X},b^{'})\to
(\overline{Y},b^{'})$. The equalities $(2.6)$ follows that the
chain maps $\overline{f}_0$ and $\overline{f}_1$ satisfy the
relations
$$\overline{f}_0(1-\overline{T}_n)=(1-\overline{T}_n)\overline{f}_1,\quad
\overline{f}_1\overline{N}_n=\overline{N}_n\overline{f}_0,\quad
n\geqslant 0.\eqno(2.7)$$ For chain bicomplexes
$(C(\overline{X}),\delta_1,\delta_2)$ and
$(C(\overline{Y}),\delta_1,\delta_2)$, consider the map of
bigraded modules $C(f):C(\overline{X})_{n,m}\to
C(\overline{Y})_{n,m}$, $n\geqslant 0$, $m\geqslant 0$, defined by
the rule $$C(f)=\left\{\begin{array}{ll}
\overline{f}_0,&\mbox{if}\quad m\equiv 0\,{\rm mod}(2),\\
\overline{f}_1,&\mbox{if}\quad m\equiv 1\,{\rm mod}(2).\\
\end{array}\right.$$
Since the chain maps $\overline{f}_0:(\overline{X},b)\to
(\overline{Y},b)$ and $\overline{f}_1:(\overline{X},b^{'})\to
(\overline{Y},b^{'})$ satisfy the relations $(2.7)$, the map of
bigraded modules $C(f):C(\overline{X})\to C(\overline{Y})$ is the
map of chain bicomlexes
$C(f):(C(\overline{X}),\delta_1,\delta_2)\to
(C(\overline{Y}),\delta_1,\delta_2)$. Now we consider the chain
complexes $(C(\overline{X}),\delta_1,\delta_2)$ and
$(C(\overline{Y}),\delta_1,\delta_2)$ together with the
above-specified action of the group
$\mathbb{Z}_2=\{1,\vartheta\,|\,\vartheta^2=1\}$. Let us show that
the map $C(f)$ is an $\mathbb{Z}_2$-equivariant map. Indeed, given
any collection $0\leqslant i_1<\dots<i_k\leqslant n$, the
conditions $(1.7)$ imply the equality
$$(-1)^{i_1+\dots+i_k}f_{(i_1,\dots,i_k)}R_n=
(-1)^{(n-i_k)+\dots+(n-i_1)}R_{n-k}f_{(n-i_k,\dots,n-i_1)}.$$
Moreover, given any collection $0\leqslant i_1<\dots<i_k\leqslant
n-1$, the relations $(1.6)$ and $(1.7)$ imply the equality
$$(-1)^{i_1+\dots+i_k}f_{(i_1,\dots,i_k)}(R_nT_n)=(-1)^{(n-i_k-1)+\dots+(n-i_1-1)}(R_{n-k}T_{n-k})
f_{(n-i_k,\dots,n-i_1)}.$$ Last two equalities follows that the
specified above families of maps $\{f_0^{\,k}\}$ and
$\{f_1^{\,k}\}$ satisfy the relations
$$f_0^{\,k}R_n=R_{n-k}f_0^{\,k},\quad
f_1^{\,k}(R_nT_n)=(R_{n-k}T_{n-k})f_1^{\,k},\quad k\geqslant
0,\quad n\geqslant 0.$$ By using these relations we obtain the
equalities
$$\overline{f}_0\overline{R}_n=\overline{R}_n\overline{f}_0,\quad
\overline{f}_1(\overline{R}_n\overline{T}_n)=(\overline{R}_n\overline{T}_n)\overline{f}_1,\quad
n\geqslant 0.$$ These equalities imply the relation
$C(f)\vartheta=\vartheta C(f)$. Thus the map of chain bicomplexes
$C(f)$ is a $\mathbb{Z}_2$-equi\-va\-riant map of chain
bicomplexes and, consequently, the map $C(f)$ induces the map of
the hyperhomology
$$H(\mathbb{Z}_2;C(f)):H(\mathbb{Z}_2;(C(\overline{X}),\delta_1,\delta_2))\to
H(\mathbb{Z}_2;(C(\overline{Y}),\delta_1,\delta_2)).$$ Thus, each
morphism of $D\F$-modules $f:(X,d,\p,t,r)\to (Y,d,\p,t,r)$ defines
the map of the dihedral homology graded modules
$$HD(f)=H(\mathbb{Z}_2;C(f)):HD(X)\to HD(Y).$$

Now, we consider the composition $gf:(X,d,\p,t,r)\to (Z,d,\p,t,r)$
of morphisms of $D\F$-modules $f:(X,d,\p,t,r)\to (Y,d,\p,t,r)$ and
$g:(Y,d,\p,t,r)\to (Z,d,\p,t,r)$. Let us show that there is the
equality of maps $HD(gf)=HD(g)HD(f)$. Consider the family of maps
$(gf)_q^k:X_{n,\bu}\to Z_{n,\bu}$, $k\geqslant 0$, $q\geqslant 0$,
defined by $(2.4)$. By using the relations $(1.5)$ we obtain that
the following relations holds: $$(gf)_q^k=\sum\limits_{0\leqslant
i_1<\dots<i_k\leqslant n-q
}\,\sum_{\sigma\in\Sigma_k}\sum_{I'_{\sigma}} (-1)^{{\rm
sign}(\sigma)+i_1+\dots+i_k}
g_{(\widehat{\sigma(i_1)},\dots,\widehat{\sigma(i_m)})}
f_{(\widehat{\sigma(i_{m+1})},\dots,\widehat{\sigma(i_k)})},$$
where $I'_{\sigma}$ is the same as in $(1.5)$. It is clear that
$$i_1+\dots+i_k+{\rm sign}(\sigma)\equiv
\widehat{\sigma(i_1)}+\dots+\widehat{\sigma(i_k)}\,{\rm mod}(2).$$
Moreover, given any collections $$0\leqslant
j_1<\dots<j_m\leqslant n-q\quad \mbox{and}\quad 0\leqslant
j_{m+1}<\dots<j_k\leqslant n-q,$$ there are a collection
$0\leqslant i_1<\dots<i_k\leqslant n-q$ and a permutation
$\sigma\in\Sigma_k$ such that $\widehat{\sigma(i_s)}=j_s$,
$1\leqslant s\leqslant k$. Therefore, the last above-specified
relations imply the relations
$$(gf)_q^k=\sum_{i+j=k}g_q^if_q^j,\quad k\geqslant 0,\quad
q\geqslant 0.$$ For $q=0,1$, by using these relations we obtain
the equalities $\overline{(gf)}_0=\overline{g}_0\overline{f}_0$
and $\overline{(gf)}_1=\overline{g}_1\overline{f}_1$. These
equalities imply the equality $C(gf)=C(g)C(f)$ of maps of chain
bicomlexes. It is clear that the map $C(gf)$ is a
$\mathbb{Z}_2$-equivariant map. Since the hyperhomology
$H(\mathbb{Z}_2;-)$ is a functor by the second argument, the
obtained above equality $C(gf)=C(g)C(f)$ implies the required
equality $HD(gf)=HD(g)HD(f)$ of maps of graded dihedral homology
modules.

Given any homotopy $h$ between morphisms $f,g:(X,d,\p,t,r)\to
(X,d,\p,t,r)$ of $D\F$-modules, in the same way as above, we
obtain the $\mathbb{Z}_2$-equi\-va\-riant homotopy
$C(h):C(\overline{X})_{*,\bu}\to C(\overline{Y})_{*+1,\bu}$
between $\mathbb{Z}_2$-equi\-va\-riant maps $C(f)$ and $C(g)$ of
chain bicomplexes. Since the functor of the hyperhomology
$H(\mathbb{Z}_2;-)$ sends homotopic $\mathbb{Z}_2$-equi\-va\-riant
maps of chain bicomplexes into equal maps of graded modules, we
obtain the equality $HD(f)=HD(g)$ of maps of graded dihedral
homology modules. This equality implies that if the map
$f:(X,d,\p,t,r)\to (X,d,\p,t,r)$ is a homotopy equivalence of
$D\F$-modules, then the induces map $HD(\widetilde{f}):HD(X)\to
HD(y)$ of graded dihedral homology modules is a isomorphism of
graded modules.~~~$\blacksquare$
\vspace{0.5cm}

\centerline{\bf \S\,3. The homotopy invariance of dihedral
homology} \centerline{\bf of involutive $\A$-algebras.}
\vspace{0.5cm}

First, following \cite{Kad} and \cite{Smir} (see also \cite{S}),
we recall necessary definitions related to the notion of an
$\A$-algebra.

An $\A$-algebra $(A,d,\pi_n)$ is any differential module $(A,d)$
with $A=\{A_n\}$, $n\in\mathbb{Z}$, $n\geqslant 0$, $d:A_\bu\to
A_{\bu-1}$, equipped with a family of maps $\{\pi_n:(A^{\otimes
(n+2)})_\bu\to A_{\bu+n}\}$, $n\geqslant 0$, satisfying the
following relations for any integer $n\geqslant -1$:
$$d(\pi_{n+1})=\sum\limits_{m=0}^n\sum_{t=1}^{m+2}(-1)^{t(n-m+1)+n+1}
\pi_m(\underbrace{1\otimes\dots\otimes
1}_{t-1}\otimes\,\pi_{n-m}\otimes\underbrace{1\otimes\dots \otimes
1}_{m-t+2}),\eqno(3.1)$$ where
$d(\pi_{n+1})=d\pi_{n+1}+(-1)^n\pi_{n+1}d$. For example, at
$n=-1,0,1$ the relations $(3.1)$ take the forms $$d(\pi_0)=0,\quad
d(\pi_1)=\pi_0(\pi_0\otimes 1)-\pi_0(1\otimes\pi_0),$$
$$d(\pi_2)=\pi_0(\pi_1\otimes
1)+\pi_0(1\otimes\pi_1)-\pi_1(\pi_0\otimes 1\otimes
1)+\pi_1(1\otimes\pi_0\otimes 1)-\pi_1(1\otimes 1\otimes\pi_0).$$

A morphism of $\A$-algebras $f:(A,d,\pi_n)\to (A',d,\pi_n)$ is
defined as a family of module maps
$f=\{f_n:(A^{\otimes(n+1)})_\bu\to
A'_{\bu+n}~|~n\in\mathbb{Z},~n\geqslant 0\}$, which, for all
integers $n\geqslant -1$, satisfy the relations
$$d(f_{n+1})=\sum_{m=0}^n\sum_{t=1}^{m+1}(-1)^{t(n-m+1)+n+1}f_m(\underbrace{1\otimes\dots\otimes
1}_{t-1}\otimes\,\pi_{n-m}\otimes\underbrace{1\otimes\dots\otimes
1}_{m-t+1})\,-$$
$$-\sum_{m=0}^n\sum_{\,\,\,J_{n-m}}(-1)^{\varepsilon}\pi_m(f_{n_1}\otimes
f_{n_2}\otimes\dots\otimes f_{n_{m+2}}),\eqno(3.2)$$ where
$d(f_{n+1})=df_{n+1}+(-1)^nf_{n+1}d$ and $$J_{n-m}=\{n_1\geqslant
0,n_2\geqslant 0,\dots,n_{m+2}\geqslant
0~|~n_1+n_2+\dots+n_{m+2}=n-m\};$$
$$\varepsilon=\sum_{i=1}^{m+1}(n_i+1)(n_{i+1}+\dots+n_{m+2}).$$
For example, at $n=-1,0,1$ the relations $(3.2)$ take,
respectively, the following view $$d(f_0)=0,\quad
d(f_1)=f_0\pi_0-\pi_0(f_0\otimes f_0),$$
$$d(f_2)=f_0\pi_1-f_1(\pi_0\otimes
1)+f_1(1\otimes\pi_0)-\pi_0(f_1\otimes f_0)+\pi_0(f_0\otimes
f_1)-\pi_1(f_0\otimes f_0\otimes f_0).$$

Under a composition of morphisms of $\A$-algebras
$f:(A,d,\pi_n)\to(A',d,\pi_n)$ and
$g:(A',d,\pi_n)\to(A'',d,\pi_n)$ we mean the morphism of
$\A$-algebras $$gf=\{(gf)_n\}:(A,d,\pi_n)\to(A'',d,\pi_n)$$
defined by
$$(gf)_{n+1}=\sum_{m=-1}^n\sum_{J_{n-m}}(-1)^\varepsilon
g_{m+1}(f_{n_1}\otimes f_{n_2}\otimes\dots\otimes
f_{n_{m+2}}),\quad n\geqslant -1,\eqno(3.3)$$ where $J_{n-m}$ and
${\varepsilon}$ are the same as in $(3.2)$. For example, at
$n=0,1,2$ the formulae $(3.3)$ take, respectively, the following
view $$(gf)_0=g_0f_0,\quad (gf)_1=g_0f_1+g_1(f_0\otimes f_0),$$
$$(gf)_2=g_0f_2-g_1(f_0\otimes f_1)+g_1(f_1\otimes
f_0)+g_2(f_0\otimes f_0\otimes f_0).$$

A homotopy between morphisms of $\A$-algebras $f,g:(A,d,\pi_n)\to
(A',d,\pi_n)$ is defined as a family of module maps
$h=\{h_n:(A^{\otimes(n+1)})_\bu\to
A'_{\bu+n+1}~|~n\in\mathbb{Z},~n\geqslant 0\}$, which, for all
integers $n\geqslant -1$, satisfy the relations
$$d(h_{n+1})=f_{n+1}-g_{n+1}+\sum_{m=0}^n\sum_{t=1}^{m+1}(-1)^{t(n-m+1)+n}h_m(\underbrace{1\otimes\dots\otimes
1}_{t-1}\otimes\,\pi_{n-m}\otimes\underbrace{1\otimes\dots\otimes
1}_{m-t+1})\,+$$
$$+\sum_{m=0}^n\sum_{\,\,\,J_{n-m}}\sum_{i=1}^{m+2}(-1)^\varrho\pi_m(g_{n_1}\otimes\dots\otimes
g_{n_{i-1}}\otimes h_{n_i}\otimes f_{n_{i+1}}\otimes\dots\otimes
f_{n_{m+2}}),\eqno(3.4)$$ where
$d(h_{n+1})=dh_{n+1}+(-1)^{n+1}h_{n+1}d$ and $J_{n-m}$ is the same
as in $(3.2)$;
$$\varrho=m+\sum_{k=1}^{m+1}(n_k+1)(n_{k+1}+\dots+n_{m+2})+\sum_{k=1}^{i-1}n_k.$$
For example, at $n=-1,0,1$ the relations $(3.4)$ take,
respectively, the following view $$d(h_0)=f_0-g_0,\quad
d(h_1)=f_1-g_1-h_0\pi_0+\pi_0(h_0\otimes f_0)+\pi_0(g_0\otimes
f_0),$$ $$d(h_2)=f_2-g_2-h_0\pi_1+h_1(\pi_0\otimes
1)-h_1(1\otimes\pi_0)-\pi_0(h_0\otimes f_1)-\pi_0(g_0\otimes
h_1)\,+$$ $$+\,\pi_0(h_1\otimes f_0)-\pi_0(g_1\otimes
h_0)-\pi_1(h_0\otimes f_0\otimes f_0)-\pi_1(g_0\otimes h_0\otimes
f_0)-\pi_1(g_0\otimes g_0\otimes h_0).$$

The origin of the signs in formulae $(3.1)$\,--\,$(3.4)$ is
described in detail in \cite{Lapin3}.

Now, we recall \cite{B} (see also \cite{SLapin1}) that an
involutive $\A$-algebra $(A,d,\pi_n,*)$ is defined as any
$\A$-algebra $(A,d,\pi_n)$ together with the automorphism of
graded modules $*:A_\bu\to A_\bu$ (the notation $*(a)=a^*$ for
$a\in A$), which, at any elements $a\in A$ and
$a_0,a_1,\dots,a_n,a_{n+1}\in A$, satisfies the conditions
$$a^{**}=a,\quad d(a^*)=d(a)^*,$$ $$\pi_n(a_0\otimes
a_1\otimes\dots\otimes a_n\otimes
a_{n+1})^*=(-1)^{\varepsilon_{n+1}}\pi_n(a_{n+1}^*\otimes
a_n^*\otimes\dots\otimes a_1^*\otimes a_0^*),\quad n\geqslant
0,\eqno(3.5)$$ where
$$\varepsilon_{n+1}=\frac{n(n+1)}{2}+\sum_{0\leqslant i<j\leqslant
n+1}|a_i||a_j|$$ and $|a|=q$ means that $a\in A_q$.

{\bf Definition 3.1}. We define a morphism of involutive
$\A$-algebras $$f:(A,d,\pi_n,*)\to (A',d,\pi_n,*)$$ as any
morphism of $\A$-algebras $f:(A,d,\pi_n)\to (A',d,\pi_n)$, which,
at any elements $a_0,a_1,\dots,a_n\in A$, satisfies the conditions
$$f_n(a_0\otimes a_1\otimes\dots\otimes a_n\otimes
a_n)^*=(-1)^{\varepsilon_n} f_n(a_n^*\otimes
a_{n-1}^*\otimes\dots\otimes a_1^*\otimes a_0^*),\quad n\geqslant
0,\eqno(3.6)$$ where
$$\varepsilon_n=\frac{(n-1)n}{2}+\sum_{0\leqslant i<j\leqslant
n}|a_i||a_j|.$$

By using the fact that any morphism of involutive $\A$-algebras is
a morphism of $\A$-algebras we define the composition of morphisms
of involutive $\A$-algebras as a composition of morphisms of
$\A$-algebras. By using the formulae $(3.3)$ and $(3.6)$ it is
easy to check that the composition of morphisms of involutive
$\A$-algebras is a morphism of involutive $\A$-algebras. Clearly,
the composition of morphisms of involutive $\A$-algebras is an
associative operation. Moreover, for each involutive $\A$-algebra
$(A,d,\pi_n,*)$, there is the identity morphism of involutive
$\A$-algebras $$1_A=\{(1_A)_n\}:(A,d,\pi_n,*)\to (A,d,\pi_n,*),$$
where $(1_A)_0={\rm id_A}$ and $(1_A)_n=0$, $n\geqslant 1$. Thus,
the class of all involutive $\A$-algebras over any commutative
unital ring $K$ and their morphisms is a category, which we denote
by $\A^{\rm{inv}}(K)$.

{\bf Definition 3.2}. We define a homotopy between morphisms of
involutive $\A$-algebras $f,g:(A,d,\pi_n,*)\to (A',d,\pi_n,*)$ as
any homotopy $h=\{h_n\}$ between morphisms of $\A$-algebras
$f,g:(A,d,\pi_n)\to (A',d,\pi_n)$, which, at any elements
$a_0,a_1,\dots,a_n\in A$, satisfies the conditions
$$h_n(a_0\otimes a_1\otimes\dots\otimes a_n\otimes
a_n)^*=(-1)^{\varepsilon_n} h_n(a_n^*\otimes
a_{n-1}^*\otimes\dots\otimes a_1^*\otimes a_0^*),\quad n\geqslant
0.\eqno(3.7)$$ where $\varepsilon_n$ is the same as in $(3.6)$.

For any involutive $\A$-algebras $(A,d,\pi_n,*)$ and
$(A',d,\pi_n,*)$, the relation between morphisms of involutive
$\A$-algebras of the form $(A,d,\pi_n,*)\to (A',d,\pi_n,*)$
defined by the presence of a homotopy between them is an
equivalence relation. By using this equivalence relation between
morphisms of involutive $\A$-algebras the notion of a homotopy
equivalence of involutive $\A$-algebras is introduced in the usual
way. Namely, a morphism of involutive $\A$-algebras is called a
homotopy equivalence of involutive $\A$-algebras, when this
morphism have a homotopy inverse morphism of involutive
$\A$-algebras.

Now, let us proceed to dihedral homology of involutive
$\A$-algebras. In \cite{SLapin1} it was shown that any involutive
$\A$-algebra $(A,d,\pi_n,*)$ defines the tensor $D\F$-module
$({\cal M}(A),d,\p,t,r)$, which is given by
$${\cal M}(A)=\{{\cal
M}(A)_{n,m}\},\quad {\cal M}(A)_{n,m}=(A^{\otimes (n+1)})_m,\quad
n\geqslant 0,\quad m\geqslant 0,$$ $$t_n(a_0\otimes\dots\otimes
a_n)=(-1)^{|a_n|(|a_0|+\dots+|a_{n-1}|)}a_n\otimes a_0\otimes
a_1\otimes \dots\otimes a_{n-1},$$ $$r_n(a_0\otimes\dots\otimes
a_n)=(-1)^{\sum_{0<i<j\leqslant n}|a_i||a_j|}a_0^*\otimes
a_n^*\otimes a_{n-1}^*\dots\otimes a_1^*,$$
$$d(a_0\otimes\dots\otimes
a_n)=\sum_{i=0}^n(-1)^{|a_0|+\dots+|a_{i-1}|}a_0\otimes\dots\otimes
a_{i-1}\otimes d(a_i)\otimes a_{i+1}\otimes\dots\otimes a_n,$$ and
the family of module maps $$\p=\{\p_{(i_1,\dots,i_k)}:{\cal
M}(A)_{n,p}\to {\cal M}(A)_{n-k,p+k-1}\},$$ $$n\geqslant 0,\quad
p\geqslant 0,\quad 1\leqslant k\leqslant n,\quad 0\leqslant
i_1<\dots<i_k\leqslant n,$$ is defined by
$$\p_{(i_1,\dots,i_k)}\!=$$ $$=\left\{\begin{array}{ll}
(-1)^{k(p-1)}1^{\otimes j}\otimes\,\pi_{k-1}\otimes 1^{\otimes
(n-k-j)}\,,\,\,\mbox{if}\,\,0\leqslant j\leqslant
n-k&\\~~~\mbox{and}\,\,\, (i_1,\dots,i_k)=(j,j+1,\dots,j+k-1);\\
(-1)^{q(k-1)}\p_{(0,1,\dots,k-1)}t_n^q\,,\,\,\mbox{if}\,\,1\leqslant
q\leqslant k&\\~~~\mbox{and}\,\,\,
(i_1,\dots,i_k)=(0,1,\dots,k-q-1,n-q+1,n-q+2,\dots,n);\\
0,\quad\mbox{otherwise}.&\\
\end{array}\right.\eqno(3.8)$$

Now we recall \cite{SLapin1} that the dihedral homology $HD(A)$ of
an involutive $\A$-algebra $(A,d,\pi_n,*)$ is defined as the
dihedral homology $HD({\cal M}(A))$ of the corresponding
$D\F$-module $({\cal M}(A),d,\p,t,r)$.

In the following assertion we describe functorial and homotopy
properties of the dihedral homology of involutive
$\A$-al\-geb\-ras.

{\bf Theorem 3.1}. The dihedral homology of involutive
$\A$-algebras over an arbitrary commutative unital ring $K$
determines the functor $HD:\A^{\rm{inv}}(K)\to GrM(K)$ from the
category of involutive $\A$-al\-gebras $\A^{\rm{inv}}(K)$ to the
category of graded $K$-modules $GrM(K)$. This functor sends
homotopy equivalences of involutive $\A$-algebras into
isomorphisms of graded modules.

{\bf Proof}. First, show that every morphism of involutive
$\A$-algebras induces a morphism of corresponding tensor
$D\F$-modules. Given any morphism of involutive $\A$-al\-geb\-ras
$f:(A,d,\pi_n,*)\to (A',d,\pi_n,*)$, we define the family of
module maps $${\cal M}(f)=\{{\cal M}(f)^n_{(i_1,\dots,i_k)}:{\cal
M}(A)_{n,p}\to {\cal M}(A')_{n-k,p+k}\},$$ $$n\geqslant 0,\quad
p\geqslant 0,\quad 0\leqslant k\leqslant n,\quad 0\leqslant
i_1<\dots<i_k\leqslant n,$$ by the following rules:

1). If $k=0$, then $${\cal M}(f)^n_{(\,\,)}=f_0^{\otimes
(n+1)}.\eqno(3.9)$$

2). If $i_k<n$ and
$(i_1,\dots,i_k)=((j^1_1,\dots,j^1_{n_1}),(j^2_1,\dots,j^2_{n_2}),\dots,(j^s_1,\dots,j^s_{n_s}))$,
$$1\leqslant s\leqslant k,\quad n_1\geqslant 1,\dots,n_s\geqslant
1,\quad n_1+\dots+n_s=k,$$ $$j_{b+1}^a=j_b^a+1,\quad 1\leqslant
a\leqslant s,\quad 1\leqslant b\leqslant n_a-1,\quad
j_1^{c+1}\geqslant j^c_{n_c}+2,\quad 1\leqslant c\leqslant s-1,$$
then $${\cal
M}(f)^n_{(i_1,\dots,i_k)}=(-1)^{k(p-1)+\gamma}\underbrace{f_0\otimes\dots\otimes
f_0}_{k_1}\otimes\,
f_{n_1}\!\otimes\underbrace{f_0\otimes\dots\otimes
f_0}_{k_2}\otimes \,f_{n_2}\,\otimes$$
$$\otimes\underbrace{f_0\otimes\dots\otimes
f_0}_{k_3}\otimes\,f_{n_3}\!\otimes\dots\otimes\underbrace{f_0\otimes\dots\otimes
f_0}_{k_s}\otimes\,f_{n_s}\!\otimes\underbrace{f_0\otimes\dots\otimes
f_0}_{k_{s+1}},\eqno(3.10)$$ where $$k_1=j_1^1,\quad
k_i=j^{\,i}_1-j^{\,i-1}_{n_{i-1}}-2,\quad 2\leqslant i\leqslant
s,\quad k_{s+1}=n+1-(k_1+\dots+k_s)-k-s,$$
$$\gamma=\sum\limits_{i=1}^{s-1}n_i(n_{i+1}+\dots+n_s).$$

3). If $i_k=n$ and
$$(i_1,\dots,i_k)=((0,1,\dots,z-1-q),(j^1_1-q,\dots,j^1_{n_1}-q),(j^2_1-q,\dots,j^2_{n_2}-q),\dots$$
$$\dots,(j^s_1-q,\dots,j^s_{n_s}-q),(n-q+1,n-q+2,\dots,n)),$$
$$z\geqslant 1,\quad  1\leqslant q\leqslant z,\quad 0\leqslant
s\leqslant k-1,\quad n_1\geqslant 1,\dots,n_s\geqslant 1,$$
$$z+n_1+\dots+n_s=k,\quad j_{b+1}^a=j_b^a+1,\quad 1\leqslant
a\leqslant s,\quad 1\leqslant b\leqslant n_a-1,$$ $$j_1^1\geqslant
z+1,\quad j_1^{c+1}\geqslant j^c_{n_c}+2,\quad 1\leqslant
c\leqslant s-1,\quad j^s_{n_s}\leqslant n-1,$$ then $${\cal
M}(f)^n_{(i_1,\dots,i_k)}=(-1)^{q(z-1)}{\cal
M}(f)^n_{((0,1,\dots,z-1),(j^1_1,\dots,j^1_{n_1}),\dots,(j^s_1,
\dots,j^s_{n_s}))}t_n^q.\eqno(3.11)$$

For example, if we consider the map ${\cal
M}(f)^{15}_{((2,3),(6,7,8))}:{\cal M}(A)_{15,p}\to {\cal
M}(A')_{10,p+5}$, then  by $(3.10)$ we obtain $${\cal
M}(f)^{15}_{((2,3),(6,7,8))}=(-1)^{5(p-1)+2\cdot 3}f_0\otimes
f_0\otimes f_2\otimes f_0\otimes
f_3\otimes\underbrace{f_0\otimes\dots\otimes f_0}_6.$$ If we
consider the map ${\cal M}(f)^n_{((0,1),(3,4),(n-2,n-1,n)}:{\cal
M}(A)_{n,p}\to {\cal M}(A')_{n-7,p+7}$, $n\geqslant 8$, then by
$(3.11)$ and $(3.10)$ we obtain $${\cal
M}(f)^n_{((0,1),(3,4),(n-2,n-1,n)}=$$ $$=(-1)^{3(5-1)}{\cal
M}(f)^n_{((0,1,2,3,4),(6,7))}t^3_n=(-1)^{7(p-1)+5\cdot
2}(f_5\otimes f_2\otimes\underbrace{f_0\otimes\dots\otimes
f_0}_{n-8})t^3_n.$$

It is worth noting that any collection of integers
$(i_1,\dots,i_k)$, $0\leqslant i_1<\dots<i_k\leqslant n$, always
can be written in the form specified in the rule 2) or in the rule
3).

Now we show that the above-specified family of maps ${\cal
M}(f)=\{{\cal M}(f)^n_{(i_1,\dots,i_k)}\}$ is the morphism of
$D\F$-modules $${\cal M}(f):({\cal M}(A),d,\p,t,r)\to ({\cal
M}(A'),d,\p,t,r).$$ In \cite{Lapin4} it was shown that defined by
the formulae $(3.9)$\,--\,$(3.11)$ the family of maps ${\cal
M}(f)=\{{\cal M}(f)^n_{(i_1,\dots,i_k)}\}$ is the morphism of
$C\F$-modules $${\cal M}(f):({\cal M}(A),d,\p,t)\to ({\cal
M}(A'),d,\p,t).$$ Therefore, we need to check that the maps ${\cal
M}(f)^n_{(i_1,\dots,i_k)}\in {\cal M}(f)$ satisfy the relations
$(1.7)$. Clearly, at $k=0$, the the equality ${\cal
M}(f)^n_{(\,\,)}r_n=r_n{\cal M}(f)^n_{(\,\,)}$ is true.  Consider
the maps ${\cal M}(f)^n_{(i_1,\dots,i_k)}$, where $i_1\geqslant 0$
and $i_k<n$, that are given by $(3.10)$ at $s=1$. Suppose that
$$(i_1,\dots,i_k)=(j^1_1,\dots,j^1_{n_1})=(j,j+1,\dots,j+k-1),\quad\mbox{where}\quad
0\leqslant j\leqslant n-k.$$ In this case, on the one hand, at any
element $a_0\otimes\dots\otimes a_n\in {\cal
M}(A)_{n,p}=(A^{n+1})_p$, we have the equalities $${\cal
M}(f)^n_{(j,j+1,\dots,j+k-1)}r_n(a_0\otimes\dots\otimes
a_n)=(-1)^\alpha {\cal M}(f)^n_{(j,j+1,\dots,j+k-1)}(a_0^*\otimes
a_n^*\otimes\dots\otimes a_1^*)=$$
$$=(-1)^{\alpha+k(p-1)}(f_{0}^{\otimes j}\otimes f_k\otimes
f_0^{\otimes (n-k-j)})(a_0^*\otimes a_n^*\otimes\dots\otimes
a_1^*)=(-1)^{\alpha+k(p-1)+\beta}f_0(a_0^*)\,\otimes$$
$$\otimes\,f_0(a_n^*)\otimes\dots\otimes f_0(a_{n-j+2}^*)\otimes
f_k(a_{n-j+1}^*\otimes\dots\otimes a_{n-j-k+1}^*)\otimes
f_0(a_{n-j-k}^*)\otimes\dots\otimes f_0(a_1^*),$$ where
$$\alpha=\sum_{0<i<j\leqslant n}|a_i||a_j|,\quad
\beta=k(|a_0|+|a_n|+|a_{n-1}|+\dots+|a_{n-j+2}|).$$ On the other
hand, we have the equalities $$(-1)^{k(k-1)/2}r_{n-k}{\cal
M}(f)^n_{(n-j-k+1,n-j-k+2,\dots,n-j)}(a_0\otimes\dots\otimes
a_n)=$$ $$=(-1)^{(k(k-1)/2)+k(p-1)}r_{n-k}(f_0^{\otimes
(n-j-k+1)}\otimes f_k\otimes f_0^{\otimes
(j-1)})(a_0\otimes\dots\otimes a_n)=$$
$$=(-1)^{(k(k-1)/2)+k(p-1)+\gamma}r_{n-k}(f_0(a_0)\otimes\dots\otimes
f_0(a_{n-j-k})\otimes f_k(a_{n-j-k+1}\otimes\dots\otimes
a_{n-j+1})\,\otimes$$ $$\otimes\,f_0(a_{n-j+2})\otimes\dots\otimes
f_0(a_n))=$$
$$=(-1)^{(k(k-1)/2)+k(p-1)+\gamma+\delta}f_0(a_0)^*\otimes
f_0(a_n)^*\otimes\dots\otimes f_0(a_{n-j+2})^*\otimes$$
$$\otimes\,f_k(a_{n-j-k+1}\otimes\dots\otimes a_{n-j+1})^*\otimes
f_0(a_{n-j-k})^*\otimes\dots\otimes f_0(a_1)^*=$$
$$=(-1)^{(k(k-1)/2)+k(p-1)+\gamma+\delta+\mu}f_0(a_0^*)\otimes
f_0(a_n^*)\otimes\dots\otimes f_0(a_{n-j+2}^*)\,\otimes$$
$$\otimes\,f_k(a_{n-j+1}^*\otimes\dots\otimes
a_{n-j-k+1}^*)\otimes f_0(a_{n-j-k}^*)\otimes\dots\otimes
f_0(a_1^*),$$ where $$\gamma=k(|a_0|+|a_1|+\dots+|a_{n-j-k}|),$$
$$\delta=\sum\limits_{0<i<j\leqslant
n}|a_i||a_j|-\sum\limits_{n-j-k+1\leqslant s<t\leqslant
n-j+1}|a_s||a_t|\,+$$
$$+\,k(|a_1|+\dots+|a_{n-j-k}|+|a_{n-j+2}|+\dots+|a_n|),$$
$$\mu=\frac{k(k-1)}{2}+\sum\limits_{n-j-k+1\leqslant s<t\leqslant
n-j+1}|a_s||a_t|.$$ Since $\alpha+\beta\equiv
(k(k-1)/2)+\gamma+\delta+\mu\,{\rm mod}(2)$, we obtain the
required relation $${\cal
M}(f)^n_{(j,j+1,\dots,j+k-1)}r_n=(-1)^{k(k-1)/2}r_{n-k}{\cal
M}(f)^n_{(n-j-k+1,n-j-k+2,\dots,n-j)}.$$ In the similar way it is
checked that relations $(1.7)$ holds for all maps ${\cal
M}(f)^n_{(i_1,\dots,i_k)}$, where $i_1\geqslant 0$ and $i_k<n$,
which are defined by $(3.10)$ at $s\geqslant 2$. Now, consider the
maps ${\cal M}(f)^n_{(i_1,\dots,i_k)}$, where $i_1\geqslant 0$ and
$i_k=n$, that are given by the formulae $(3.11)$ at $s=0$. Suppose
that
$$(i_1,\dots,i_k)=((0,1,\dots,k-1-q),(n-q+1,n-q+2,\dots,n)),$$
where $1\leqslant q\leqslant k$. In this case by $(3.11)$, at
$s=0$, we have $${\cal
M}(f)^n_{(0,1,\dots,k-q-1,n-q+1,n-q+2,\dots,n)}=(-1)^{q(k-1)}{\cal
M}(f)^n_{(0,1,\dots,k-1)}t_n^q.$$ By using the relations
$t_n^qr_n=r_nt_n^{-q}=r_nt_n^{n+1-q}=r_nt_n^{n-k+1}t_n^{k-q}$ we
obtain $${\cal
M}(f)^n_{(0,1,\dots,k-q-1,n-q+1,n-q+2,\dots,n)}r_n=(-1)^{q(k-1)}{\cal
M}(f)^n_{(0,1,\dots,k-1)}t_n^qr_n=$$ $$=(-1)^{q(k-1)}{\cal
M}(f)^n_{(0,1,\dots,k-1)}r_nt_n^{n-k+1}t_n^{k-q}=$$ $$=
(-1)^{q(k-1)+(k(k-1)/2)}r_{n-k}{\cal
M}(f)^n_{(n-k+1,\dots,n)}t_n^{n-k+1}t_n^{k-q}.$$ In \cite{Lapin4}
it was shown that the maps ${\cal M}(f)^n_{(i_1,\dots,i_k)}\in
{\cal M}(f)$ satisfy the relations $(1.6)$. Therefore, we obtain
$$(-1)^{q(k-1)+(k(k-1)/2)}r_{n-k}{\cal
M}(f)^n_{(n-k+1,\dots,n)}t_n^{n-k+1}t_n^{k-q}=$$
$$=(-1)^{q(k-1)+(k(k-1)/2)}r_{n-k}t_{n-k}^{n-k+1}{\cal
M}(f)^n_{(0,1,\dots,k-1)}t_n^{k-q}=$$
$$=(-1)^{q(k-1)+(k(k-1)/2)}r_{n-k}{\cal
M}(f)^n_{(0,1,\dots,k-1)}t_n^{k-q}.$$ By using the formulae
$(3.11)$ we have $${\cal
M}(f)^n_{(0,1,\dots,q-1,n-k+q+1,\dots,n)}=(-1)^{(k-q)(k-1)}{\cal
M}(f)^n_{(0,1,\dots,k-1)}t_n^{k-q}.$$ Since
$q(k-1)+(k(k-1)/2)+(k-q)(k-1)\equiv (k(k-1)/2)\,{\rm mod}(2)$, we
obtain the required relation $${\cal
M}(f)^n_{(0,1,\dots,k-q-1,n-q+1,n-q+2,\dots,n)}r_n=
(-1)^{k(k-1)/2}r_{n-k}{\cal
M}(f)^n_{(0,1,\dots,q-1,n-k+q+1,\dots,n)}.$$ In the similar way it
is checked that the relations $(1.7)$ holds for all maps ${\cal
M}(f)^n_{(i_1,\dots,i_k)}$, where $i_1\geqslant 0$ and $i_k=n$,
which are given by $(3.11)$ at $s\geqslant 1$. Thus, all module
maps ${\cal M}(f)^n_{(i_1,\dots,i_k)}\in{\cal M}(f)$ satisfy the
relations $(1.7)$. It follows that the specified above morphism of
$C\F$-modules ${\cal M}(f):({\cal M}(A),d,\p,t)\to ({\cal
M}(A'),d,\p,t)$ is the morphism of $D\F$-modules ${\cal
M}(f):({\cal M}(A),d,\p,t,r)\to ({\cal M}(A'),d,\p,t,r)$.

Now, we consider the composition $gf:(A,d,\pi_n,*)\to
(A'',d,\pi_n,*)$ of an arbitrary morphisms $f:(A,d,\pi_n,*)\to
(A',d,\pi_n,*)$ and $g:(A',d,\pi_n,*)\to (A'',d,\pi_n,*)$ of
involutive $\A$-algebras. In the same way as in \cite{Lapin4} it
is checked that the equality of morphisms of $D\F$-modules ${\cal
M}(gf)={\cal M}(g){\cal M}(f)$ is true.

The above follows that there is the functor ${\cal
M}:\A^{\rm{inv}}(K)\to D\F(K)$. The required functor
$HD:\A^{\rm{inv}}(K)\to GrM(K)$ we define as a composition of the
functor ${\cal M}:\A^{\rm{inv}}(K)\to D\F(K)$ and the functor
$HD:D\F(K)\to GrM(K)$ specified in Theorem 2.1.

Now we show that the functor $HD:\A^{\rm{inv}}(K)\to GrM(K)$ sends
homotopy equivalences of involutive $\A$-algebras into
isomorphisms of graded modules. Taking into account Theorem 2.1,
it suffices to show that the functor ${\cal M}:\A^{\rm{inv}}(K)\to
D\F(K)$ sends homotopy equivalences into homotopy equivalences of
$D\F$-modules. With this purpose we show that each homotopy
between morphisms of involutive $\A$-algebras induces a homotopy
between corresponding morphisms of tensor $D\F$-modules.

Given any homotopy $h=\{h_n:(A^{\otimes(n+1)})_\bu\to
A'_{\bu+n+1}~|~n\in\mathbb{Z},~n\geqslant 0\}$ between morphisms
of involutive $\A$-algebras $f,g:(A,d,\pi_n,*)\to (A',d,\pi_n,*)$,
we define a family of module maps $${\cal M}(h)=\{{\cal
M}(h)^n_{(i_1,\dots,i_k)}:{\cal M}(A)_{n,p}\to {\cal
M}(A')_{n-k,p+k+1}\},$$ $$n\geqslant 0,\quad p\geqslant 0,\quad
0\leqslant k\leqslant n,\quad 0\leqslant i_1<\dots<i_k\leqslant
n,$$ by the following rules:

1$'$). If $k=0$, then $${\cal
M}(h)^n_{(\,\,)}=\sum\limits_{i=1}^{n+1}\underbrace{g_0\otimes\dots\otimes
g_0}_{i-1}\otimes\,h_0\otimes\underbrace{f_0\otimes\dots\otimes
f_0}_{n-i+1}\,;$$

2$'$). If $i_k<n$ and the collection
$$(i_1,\dots,i_k)=((j^1_1,\dots,j^1_{n_1}),(j^2_1,\dots,j^2_{n_2}),\dots,(j^s_1,\dots,j^s_{n_s}))$$
is the same as in the above rule 2) defining the formula $(3.10)$,
then $${\cal
M}(h)^n_{(i_1,\dots,i_k)}=(-1)^{k(p-1)+\gamma}\sum_{i=1}^s(-1)^{n_1+\dots+n_{i-1}}\underbrace{g_0\otimes\dots\otimes
g_0}_{k_1}\otimes\,g_{n_1}\!\otimes\dots$$
$$\dots\,\otimes\underbrace{g_0\otimes\dots\otimes
g_0}_{k_{i-1}}\otimes\,g_{n_{i-1}}\!\otimes\underbrace{g_0\otimes\dots\otimes
g_0}_{k_i}\otimes\,h_{n_i}\!\otimes\underbrace{f_0\otimes\dots\otimes
f_0}_{k_{i+1}}\otimes\,f_{n_{i+1}}\!\otimes\dots$$
$$\dots\otimes\underbrace{f_0\otimes\dots\otimes
f_0}_{k_s}\otimes\,f_{n_s}\!\otimes\underbrace{f_0\otimes\dots\otimes
f_0}_{k_{s+1}}\,+$$
$$+\,(-1)^{k(p-1)+\gamma}\sum_{i=1}^{s+1}(-1)^{n_1+\dots+n_{i-1}}\underbrace{g_0\otimes\dots\otimes
g_0}_{k_1}\otimes\,g_{n_1}\!\otimes\dots\otimes\underbrace{g_0\otimes\dots\otimes
g_0}_{k_{i-1}}\otimes\,g_{n_{i-1}}\otimes$$
$$\otimes\left\{\sum_{j=1}^{k_i}\underbrace{g_0\otimes\dots\otimes
g_0}_{j-1}\otimes\,h_0\otimes\underbrace{f_0\otimes\dots\otimes
f_0}_{k_i-j}\right\}\otimes\,f_{n_i}\!\otimes\underbrace{f_0\otimes\dots\otimes
f_0}_{k_{i+1}}\otimes\,f_{n_{i+1}}\otimes\dots$$
$$\dots\otimes\underbrace{f_0\otimes\dots\otimes
f_0}_{k_s}\otimes\,f_{n_s}\!\otimes\underbrace{f_0\otimes\dots\otimes
f_0}_{k_{s+1}}\,,$$ where $k_1,\dots,k_{s+1}$ and $\gamma$ are the
same as in $(3.10)$;

3$'$). If $i_k=n$ and the collection
$$(i_1,\dots,i_k)=((0,1,\dots,z-1-q),(j^1_1-q,\dots,j^1_{n_1}-q),(j^2_1-q,\dots,j^2_{n_2}-q),\dots$$
$$\dots,(j^s_1-q,\dots,j^s_{n_s}-q),(n-q+1,n-q+2,\dots,n))$$ is
the same as in the above rule 3) defining the formula $(3.11)$,
then $${\cal M}(h)^n_{(i_1,\dots,i_k)}=(-1)^{q(z-1)}{\cal
M}(h)^n_{((0,1,\dots,z-1),(j^1_1,\dots,j^1_{n_1}),\dots,(j^s_1,
\dots,j^s_{n_s}))}t_n^q.$$

In similar way as it was done above in the case of morphisms of
$D\F$-mo\-du\-les ${\cal M}(f)$, it is proved that defined by any
homotopy $h=\{h_n\}$ between morphisms of involutive $\A$-algebras
$f,g:(A,d,\pi_n,*)\to (A',d,\pi_n,*)$ the family of module maps
${\cal M}(h)=\{{\cal M}(h)^n_{(i_1,\dots,i_k)}\}$ is a homotopy
between corresponding morphisms of $D\F$-modules  ${\cal
M}(f),{\cal M}(g):({\cal M}(A),d,\p,t,r)\to ({\cal
M}(A'),d,\p,t,r)$. It follows that if $f:(A,d,\pi_n,*)\to
(A',d,\pi_n,*)$ is a homotopy equivalence of involutive
$\A$-algebras, then the corresponding morphism ${\cal M}(f):({\cal
M}(A),d,\p,t,r)\to ({\cal M}(A'),d,\p,t,r)$ is a homotopy
equivalence of $D\F$-modules. Thus, the functor ${\cal M}:\A^{{\rm
inv}}(K)\to D\F(K)$ sends homotopy equivalences of involutive
$\A$-algebras into homotopy equivalences of $D\F$-modules and,
consequently, the functor $HD:\A^{{\rm inv}}(K)\to GrM(K)$ sends
homotopy equivalences of involutive $\A$-algebras into
isomorphisms of graded modules.~~~$\blacksquare$

\vspace{1cm}

Serov Str., Saransk, Russia,

e-mail: slapin@mail.ru

\end{document}